\documentclass[11pt]{article}

\usepackage{amsmath}
\usepackage{graphicx}
\usepackage{epsfig}
\usepackage{amsfonts}
\usepackage{amssymb}
\usepackage{placeins}
\usepackage{cases}
\usepackage[margin=1in]{geometry}
\usepackage{authblk}
\usepackage[titletoc,toc,title]{appendix}
\usepackage{epstopdf}
\usepackage{xcolor}
\usepackage[final]{changes}
\usepackage{marvosym}
\usepackage{enumerate} 

\usepackage{accents}

\title{Data-Driven Inference of High-Accuracy Isostable-Based Dynamical Models in Response to External Inputs}   

\begin{document}
\author[1]{Dan Wilson}
\affil[1]{Department of Electrical Engineering and Computer Science, University of Tennessee, Knoxville, TN 37996, USA}
\maketitle

\begin{abstract}
Isostable reduction is a powerful technique that can be used to characterize behaviors of nonlinear dynamical systems in a basis of slowly decaying eigenfunctions of the Koopman operator.  When the underlying dynamical equations are known, previously developed numerical techniques allow for high-order accuracy computation of isostable reduced models.  However, in situations where the dynamical equations are unknown, few general techniques are available that provide reliable estimates of the isostable reduced equations, especially in applications where large magnitude inputs are considered.  In this work, a purely data-driven inference strategy yielding high-accuracy isostable reduced models is developed for dynamical systems with a fixed point attractor. By analyzing steady state outputs of nonlinear systems in response to sinusoidal forcing, both isostable response functions and isostable-to-output relationships can be estimated to arbitrary accuracy in an expansion performed in the isostable coordinates.    Detailed examples are considered for a population of synaptically coupled neurons and for the one-dimensional Burgers' equation.  While linear estimates of the isostable response functions are sufficient to characterize the dynamical behavior when small magnitude inputs are considered, the high-accuracy reduced order model inference strategy proposed here is essential when considering large magnitude inputs.
\end{abstract}

\pagebreak

%
\textbf{ Data-driven identification of low-order dynamical models is of primary importance in the physical, chemical, and biological sciences.  Many well-established model identification strategies have been developed to identify linear models from data, however, such strategies often fail in situations where nonlinear terms dominate.   In this work, the Koopman operator paradigm is leveraged to develop a high-accuracy model inference strategy for dynamical systems with stable fixed points.  By analyzing steady state outputs of nonlinear systems in response to sinusoidal forcing it is shown that a general, nonlinear dynamical model can be inferred.  Input-to-output relationships of the resulting model can be obtained to arbitrary accuracy in an asymptotic expansion performed in a reduced coordinate basis.  In detailed numerical examples with applications to both neuroscience and fluid dynamics, the proposed strategy identifies a low-order model that accurately captures fundamentally nonlinear system behaviors.
}
\\

\section{Introduction}
With an ever-growing abundance of data, there is an increasing demand for data-driven strategies that can be used to both predict and control the behaviors of dynamical systems.    While a variety of projection based techniques \cite{benn15} are well-suited for dynamical systems that can be accurately approximated by a local linearization, Koopman analysis has emerged as an essential technique in the study of dynamical systems with fundamentally nonlinear behaviors.  The key feature of  Koopman analysis is that it allows the observables of a nonlinear dynamical system to be studied in terms of a linear, but possibly infinite-dimensional operator \cite{budi12}, \cite{mezi13}.  As such, complicated nonlinear behaviors can be understood in terms of the  evolution of eigenfunctions that propagate in accordance with their associated eigenvalues \cite{maur16}, \cite{bagh13}, \cite{arba17}.  Using the Koopman operator framework, it is generally possible to make predictions about the evolution of observables by considering their behavior on a low-dimensional manifold comprised of a finite number of Koopman eigenmodes \cite{lusc18}, \cite{brun17}.

While  Koopman eigenmodes can often be identified if the underlying dynamical equations are known, however, when the underlying dynamical equations are unknown, the Koopman operator must be approximated in a data-driven setting.  Dynamical mode decomposition (DMD) is one such data-driven method \cite{rowl09}, \cite{kutz16}, which approximates eigenmodes from a timeseries of system observables.  Extended DMD \cite{will15},  deep learning approaches \cite{yeun19}, \cite{lusc18}, and delay embedding \cite{arba17}, \cite{brun17}, \cite{wils20ddred} techniques have also been useful in certain applications.  A commonality among these data-driven Koopman analysis strategies is that they attempt to identify a low-dimensional representation of the infinite-dimensional Koopman operator.   

Because the Koopman framework often results in a low-dimensional dynamical representation for the evolution of system observables, there has been a growing interest in applying this strategy to control problems.  Time-varying inputs make the system behavior non-autonomous, and in such systems, the Koopman eigenmodes and associated eigenvalues are time-dependent adding an additional layer of complexity to the analysis.  With this in mind, a number of Koopman-based control strategies have been developed that augment the system measurements with knowledge of the applied control in order to characterize the dynamical behavior in response to inputs \cite{proc16}, \cite{proc18}, \cite{will16}, \cite{kord18}.  Related strategies have been successful that consider a finite number of constant controls, computing the resulting Koopman eigenfunctions, and switching between control states as necessary \cite{peit19}.  An alternative strategy is to identify Koopman eigenfunctions of the autonomous, unforced system and work directly in the resulting eigenfunction basis \cite{kais17}.  This is the approach taken by the isostable reduction framework  \cite{maur13}, \cite{wils15}, \cite{wils16isopde}, \cite{soot17}, \cite{wils17isored} that uses a basis of the most slowly decaying eigenfunctions associated with the Koopman operator.  In many applications an accurate reduced order model can be obtained using a small number of isostable coordinates \cite{wils16isopde}, \cite{wils17isored}, \cite{wils19prl} making control and analysis possible in a substantially reduced order setting.  

Characterizing the response of isostable coordinates to applied input is a fundamental challenge when considering control applications in the isostable reduced coordinate basis.  This can be done, for instance, by approximating the gradient of the isostable coordinate field with respect to the applied inputs.  Previous work has focused on identifying linear approximations of the isostable response functions with respect to a fixed point attractor \cite{wils19cdc}, \cite{wils20acc}, or along specific trajectories \cite{wils15}, \cite{wils16isopde}.  These strategies work in the limit that the applied inputs are small in magnitude but suffer when larger magnitude inputs are considered.  This limitation can be addressed by considering high order expansions of the isostable response functions in a basis of isostable coordinates.  Recent work \cite{wils20highacc} details such a method for approximating these isostable response functions in oscillatory dynamical systems  to arbitrary accuracy, however, this strategy requires knowledge of the underlying model equations.  Currently, no existing strategies have been developed for identifying high-accuracy approximations of the isostable response functions using data-driven techniques.  

In this work, a purely data-driven strategy is developed for identifying reduced order isostable-based models in dynamical systems with stable fixed points.  By analyzing  steady-state model outputs in response to sinusoidal forcing, both the isostable response functions and the isostable-to-output relationships can be estimated to arbitrary accuracy in an expansion performed in the basis of isostable coordinates.  This proposed method can be considered an extension of the technique suggested in \cite{wils19cdc}, which considered only the first order accurate dynamics of the isostable reduced equations.  The current work allows allows for estimates of higher order accuracy terms which are critically important when considering large magnitude inputs.  The organization of this paper is as follows:~Section \ref{isoback} gives necessary background on the isostable reduced coordinate framework and also reframes the numerical isostable reduction techniques from \cite{wils20highacc} for use in  dynamical systems with stable fixed points.  Section \ref{fittingsection} details a data-driven strategy for inferring isostable-based reduced order models using information about steady state model outputs in response to sinusoidal forcing.   Section \ref{exampsec} provides examples that include systems with relevance to both neuroscientific and fluid flow applications, and Section \ref{concsec} gives concluding remarks.


\section{Background:~Isostable Reduction} \label{isoback}
Consider an ordinary differential equation
\begin{equation} \label{maineq}
\dot{x} = F(x) + U(t),
\end{equation}
where $x \in \mathbb{R}^N$ denotes the state, $F$ represents the nominal dynamics, and $U(t)$ is an external input.  Let $\phi(t,x)$ denote the unperturbed flow of \eqref{maineq}, and suppose  that \eqref{maineq} has a stable fixed point denoted by $x_0$.  Near the fixed point, local linearization yields
\begin{equation} \label{locallin}
\Delta x = J \Delta x + O(|\Delta x|^2),
\end{equation}
where $\Delta x \equiv  x-x_0$ and $J$ is the Jacobian of $F$ evaluated at $x_0$.  While \eqref{locallin} is only accurate in a close vicinity of the fixed point, it nonetheless can be used in conjunction with the isostable coordinate framework to characterize the infinite-time decay of solutions in the fully nonlinear basin of attraction of $x_0$.  To do so, let $w_k$, $v_k$, and $\lambda_k$ to be left eigenvectors, right eigenvectors, and eigenvalues of $J$, respectively, ordered so that $| {\rm Re}(\lambda_k)| \leq | {\rm Re}(\lambda_{k+1})|$.   Isostable coordinates associated with a subset eigenvalues with the smallest magnitude real components can be defined explicitly according to 
\begin{equation} \label{isodef}
\psi_j(x) = \lim_{t \rightarrow \infty} \left(  w_j^T  (\phi(t,x)-x_0) \exp(-\lambda_j t) \right),
\end{equation}
where $^T$ indicates the matrix transpose.  Intuitively, in Equation \eqref{isodef}, the term   $w_j^T  (\phi(t,x)-x_0)$ selects for the component of the decay in the $v_j$ direction -- multiplying by the exponentially growing term an taking the limit as time approaches infinity yields the isostable coordinate.  The definition \eqref{isodef} is closely aligned with the one given in \cite{wils16isopde}, however, other definitions that compute Fourier averages of observables \cite{maur13}, \cite{wils15} can also be used.  As a point of emphasis, the constructive definition of isostable coordinates  \eqref{isodef} is only possible for a subset of eigenvalues with the smallest magnitude real components \cite{kval19}.   Isostable coordinates associated with larger magnitude eigenvalues must be defined implicitly as level sets of Koopman eigenfunctions with decay rates that are governed by their associated $\lambda_j$.   

 The utility of the transformation to isostable coordinates relies on the fact that many dynamical systems have behavior that is governed by their slowest decaying modes of the Koopman operator.  By truncating the rapidly decaying components and retaining only a small number of the slowly decaying isostables, a reduced order model can be obtained.  Such isostable coordinate systems have been applied recently to a variety of control applications \cite{soot17}, \cite{maur18}, \cite{cast20}, \cite{shir17}, \cite{mong19}, \cite{wils19complex}.

\subsection{Reduced Order Models Using Isostable Coordinates} \label{reducedexamp}
When $U(t) = 0$, under the evolution of $\dot{x} = F(x)$,  all isostable coordinates decay exponentially according to $\dot{\psi}_k = \lambda_k \psi_j$ in the entire nonlinear basin of attraction of the fixed point.       In the analysis to follow, the dynamics of the $M$ slowest decaying isostable coordinates will be explicitly considered and  the rest will be truncated (i.e., by taking  all truncated isostable coordinates to be equal to zero at all times). The reduced order  dynamics of Equation \eqref{maineq} on a hypersurface determined by the non-truncated isostable coordinates then follows
\begin{align} \label{reducedmod}
\dot{\psi}_k &= \lambda_k \psi_k + I_k(\psi_1,\dots,\psi_M) \cdot U(t), \nonumber \\
&k = 1,\dots,M, \nonumber \\
x  &= x_0 + G(\psi_1,\dots,\psi_M),
\end{align}
where $I_k$ denotes the gradient of $\psi_k$, $G$ maps the isostable coordinates to the state, and the dot denotes the dot product.  For low dimensional systems (with $N \leq 3$), it can sometimes be possible to compute the isostable coordinates directly for all states $x$ and infer the functions $I_k$ and $G$ numerically \cite{maur13},\cite{wils19phase}.  This approach, however, does not scale well to in higher dimensions.  Instead, it is usually more feasible to Taylor expand both $I_k$ and $G$ in powers of the isostable coordinates as
\begin{align} 
G(\psi_1,\dots,\psi_M)  &\approx \sum_{k = 1}^{M} \left[ \psi_k g^k \right] + \sum_{j = 1}^M \sum_{k = 1}^j \left[ \psi_j \psi_k g^{jk} \right]+ \sum_{i = 1}^M \sum_{j = 1}^i \sum_{k = 1}^j \left[ \psi_i \psi_j \psi_k g^{ijk}\right] + \dots, \label{taylG} \\
{I}_n(\psi_1,\dots,\psi_M) &\approx I_n^0 + \sum_{k = 1}^M \left[ \psi_k I_n^k\right]  + \sum_{j = 1}^M \sum_{k = 1}^j \left[  \psi_j \psi_k I_n^{jk}  \right] + \sum_{i = 1}^M \sum_{j = 1}^i \sum_{k = 1}^j \left[ \psi_i \psi_j \psi_k I_n^{ijk}\right] + \dots.  \label{taylI}
\end{align}
For the purposes of this work Equation \eqref{reducedmod} will be considered a $j^{\rm th}$ order accurate model if the expansion of $G$ from \eqref{taylG} and and each $I_n$ from \eqref{taylI} are computed to $j^{\rm th}$ and $j-1^{\rm th}$ order accuracy in the isostable coordinates, respectively.  This convention generally assumes that $U(t)$ and each $\psi_n$ are $O(\epsilon)$ terms where $0 < \epsilon \ll 1$ so that the resulting reduction \eqref{reducedmod} retains terms up to and including $O(\epsilon^j)$ accuracy.  Note that these models often perform well even when the when the magnitude of input is larger than $O(\epsilon)$ (see, for example, \cite{wils20highacc}).

\begin{figure}[htb]
\begin{center}
\includegraphics[height=4in]{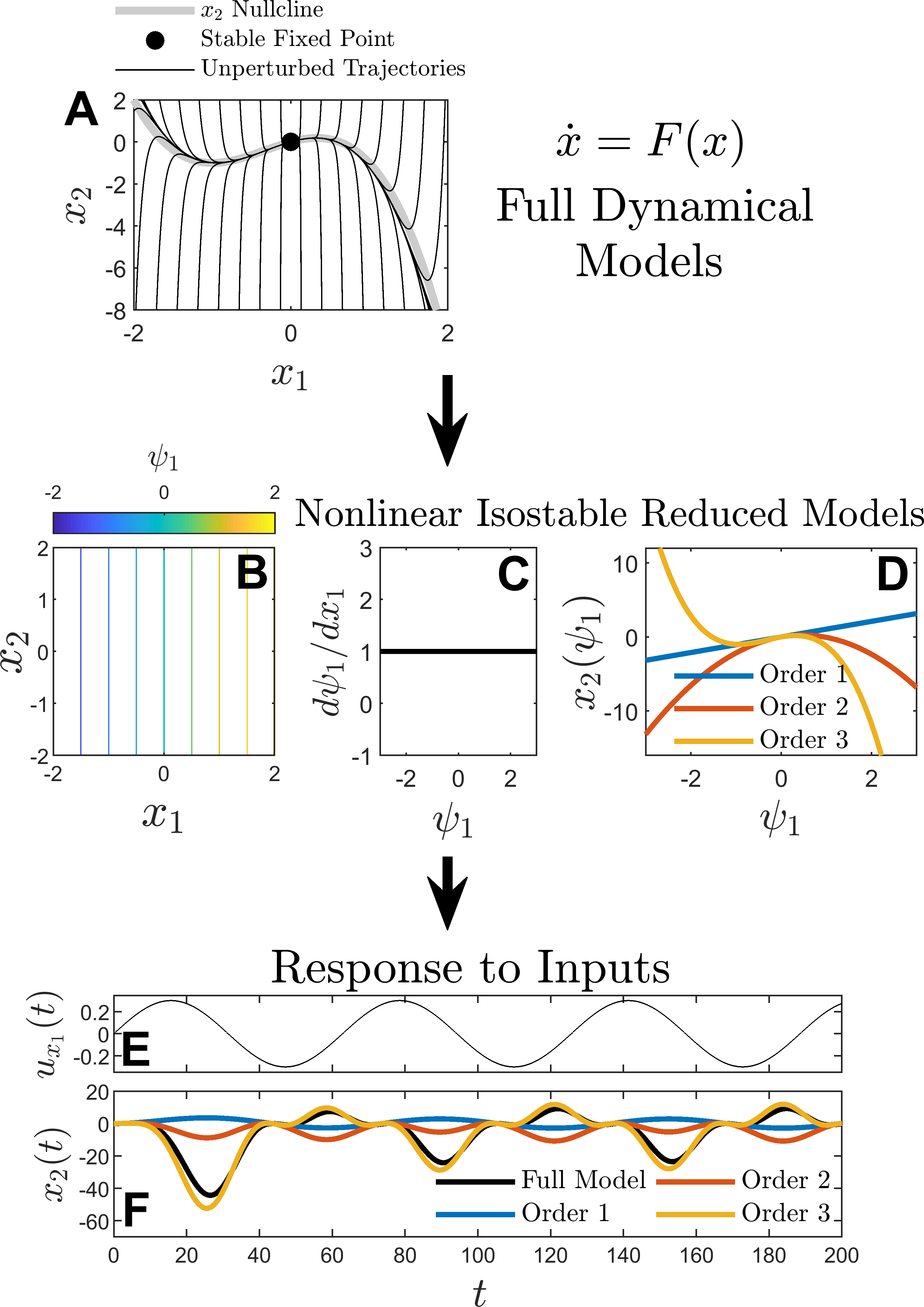}
\end{center}
\caption{A Model dimension reduction applied to a simple model \eqref{exampfn} with nominal dynamics represented in panel A.  If the full model equations are known, it is straightforward to compute the necessary terms of the reduced order models from Equations \eqref{taylG} and \eqref{taylI} that only include the slowest decaying isostable dynamics; these terms are shown in panels B through D for a single isostable coordinate.  Once the terms of the reduced order model are identified, it can be used to predict the response to inputs, provided the terms of $G$ and $I_1$ are computed to sufficiently high order accuracy in the isostable coordinate.}
\label{exampfig}
\end{figure}

Using methods adapted from \cite{wils20highacc} that are given in Appendix \ref{apxa}, it is possible to directly solve for each $g^{ijk\dots}$ and $I_n^{ijk\dots}$ directly when the underlying model equations are known.  As a concrete example, a simple dynamical system is considered
\begin{align} \label{exampfn}
\dot{x}_1 & =  \mu x_1 + u_{x_1}(t), \nonumber \\
\dot{x}_2 & = \lambda(-x_1 + x_2 + x_1^2 + x_1^3),
\end{align}
with $\lambda = -1$ and $\mu = -.05$.  The model \eqref{exampfn} has a finite-dimensional Koopman invariant subspace \cite{brun16} making it an ideal testbed for developing nonlinear model identification strategies.  Equation \eqref{exampfn} has a stable fixed point at the $x_1 = x_2 = 0$ with associated eigenvalues $\lambda_1 = -0.05$ and $\lambda_2 = -1$.  Trajectories tend to approach along the $x_2$ nullcline as shown in Figure \ref{exampfig}.  Due to the separation between the eigenvalues, a single isostable reduced model can be determined by identifying the unknown $g^{ijk}$ and $I_1^{ijk}$ terms from equations \eqref{taylG} and \eqref{taylI}, respectively, using methods detailed in Appendix \ref{apxa}.   Level sets of the isostable coordinates and the gradient of the isostable coordinates are shown in panels B and C of Figure \ref{exampfig}.  Approximations of $G(\psi_1)$ valid to various orders of $\psi_1$ in the expansion from \eqref{taylG} are shown in panel D.  The resulting reduced dimension model of the form \eqref{reducedmod} simulated in response to a sinusoidal input shown in panels E, corresponding outputs are given in panel F.   For \eqref{exampfn}, third order accuracy is sufficient to match the full model behavior.


The procedure described in Figure \ref{exampfig} explicitly requires knowledge of the underlying model equations.  In practical applications, underlying model equations are generally not known {\it a priori}.  Rather, time series of model observables are typically all that can be measured.   With this in mind, the following sections will illustrate a purely data-driven nonlinear model identification strategy.

\section{Nonlinear Model Identification from Input-Output Relationships Using Isostable Reduced Models} \label{fittingsection}
In the analysis to follow, input-output relationships will be used to identify the necessary terms of an isostable reduced model based on \eqref{reducedmod}.  To proceed, it will be assumed that $U(t) =  A u(t)$, with $A = \begin{bmatrix} A_1 & \dots & A_N  \end{bmatrix}^T$ and $u(t) \in \mathbb{R}$; in other words, $U(t)$ is explicitly assumed to be a rank-1 input.  Higher rank inputs can be considered with appropriate modifications.  Sinusoidal inputs of the form  $u(t) =  \epsilon \sin(\omega t)$ will be  used to measure the input-output dynamics  where $\omega$ gives the natural frequency of the input and $\epsilon$ is assumed to be a small constant.  Intuitively, when a sinusoidal input is applied to a {\it linear} model, the steady state output is a phase-shifted sinusoid at the same frequency.  As will be shown in the analysis to follow, a nonlinear reduced order model can be directly related to higher-frequency Fourier modes comprising the steady-state output -- this information can then be used to infer nonlinear isostable reduced models that characterize the response to inputs.

To proceed, using the assumed structure on $U(t)$, the isostable coordinates reduction from \eqref{reducedmod} will be written as
\begin{align} \label{reducedmodinput}
\dot{\psi}_k &= \lambda_k \psi_j +  I_{k,A}(\psi_1,\dots,\psi_M) \cdot u(t),  \\
&k = 1,\dots,M,  \nonumber \\
y &= y_0 + G_y(\psi_1,\dots,\psi_M)  \label{reducedoutput}.
\end{align}
Above, $y \in \mathbb{R}^\zeta$ is a collection of observables and $G_y$ maps the isostable coordinates to the system observables.  As with \eqref{reducedmod}, $G_y$ will be Taylor expanded in the basis of isostable coordinates 
\begin{equation}\label{taylgobs}
G_y(\psi_1,\dots,\psi_M)  \approx \sum_{k = 1}^{M} \left[ \psi_k g_y^k \right] + \sum_{j = 1}^M \sum_{k = 1}^j \left[ \psi_j \psi_k g_y^{jk} \right]+ \sum_{i = 1}^M \sum_{j = 1}^i \sum_{k = 1}^j \left[ \psi_i \psi_j \psi_k g_y^{ijk}\right] + \dots.
\end{equation} 
Similarly, the isostable response to inputs from \eqref{reducedmodinput} will  be written in terms of the Taylor expansion
\begin{equation} \label{tinput}
{I}_{n,A}(\psi_1,\dots,\psi_M) \approx I_{n,A}^0 + \sum_{k = 1}^M \left[ \psi_k I_{n,A}^k\right]  + \sum_{j = 1}^M \sum_{k = 1}^j \left[  \psi_j \psi_k I_{n,A}^{jk}  \right] + \sum_{i = 1}^M \sum_{j = 1}^i \sum_{k = 1}^j \left[ \psi_i \psi_j \psi_k I_{n,A}^{ijk}\right] + \dots,
\end{equation}
where, for instance, $I_{n,A}^{ijk} = I_{n}^{ijk} \cdot A$.   As with Equation \eqref{reducedmod}, the reduced order model \eqref{reducedmodinput} with output \eqref{reducedoutput} will be considered a $j^{\rm th}$ order accurate model if the terms of $G_y$ and the terms of each $I_{n,A}$ are computed to $j^{\rm th}$ and $j-1^{\rm th}$ order accuracy in the isostable coordinates, respectively.


In order to aid in the analysis,  each $\psi_n$ will be asymptotically expanded in orders of $\epsilon$ as 
\begin{equation} \label{psiexp}
\psi_n(t) = \psi_n^{(0)}(t) + \epsilon \psi_n^{(1)}(t)  + \epsilon^2 \psi_n^{(2)}(t) + \dots.
\end{equation}
Substituting \eqref{psiexp} into \eqref{reducedmodinput} and collecting the  $O(1)$ terms yields $\dot{\psi}_n^{(0)} = \lambda_n  \psi_n^{(0)}$; consequently in the limit that time approaches infinity,  $\psi_n^{(0)}$ tends to zero.  With this in mind, for simplicity of exposition, it will be assumed that each $\psi_n^{(0)}$ is zero for all time.

\subsection{Inferring First Order Accurate Models for Single Output Systems} \label{firstsec}
 For the moment, it will be assumed that $\zeta = 1$, i.e.,~that the observable $y$ is of dimension 1.  Higher dimensional outputs will be considered in Section \ref{multisec}.  A similar computation strategy for first order accurate models was considered in \cite{wils19cdc}.  To begin, $u(t) =  \epsilon \sin(\omega t)$ will be applied  and only the steady state dynamics will be considered.  Collecting the $O(\epsilon)$ terms after substituting    \eqref{psiexp} into \eqref{reducedmodinput} yields
\begin{equation}
\dot{\psi}_n^{(1)} = \lambda_n \psi_n^{(1)} + I_{n,A}^0 \sin(\omega t).  
\end{equation}
The steady state, the solution to \eqref{reducedmodinput} is 
\begin{align} \label{ord1sol}
\psi_{n,ss}^{(1)} &= \frac{- I_{n,A}^0}{ \omega^2 + \lambda_n^2}  (\lambda_n \sin(\omega t) + \omega \cos(\omega t)) \nonumber \\ 
&=   s_{1,n}^0(\omega)  I_{n,A}^0   \sin(\omega t) +   c_{1,n}^0(\omega)  I_{n,A}^0  \cos(\omega t),
\end{align}
where $s_{1,n}^0(\omega)  \in \mathbb{C}$ and $c_{1,n}^0(\omega)  \in \mathbb{C}$ are defined appropriately.  Considering \eqref{reducedoutput}, to leading order in $\epsilon$ the output is 
\begin{equation} \label{xss}
y_{ss}(w,t)  =y_0 +  \epsilon  \sum_{k = 1}^M \bigg [ s_{1,n}^0(\omega) g_y^k I_{k,A}^0   \sin(\omega t) +   c_{1,n}^0(\omega) g_y^k I_{k,A}^0  \cos(\omega t) \bigg ]+ O(\epsilon ^2),
\end{equation}
Multiplying both sides of \eqref{xss} by either $\sin(\omega t)$ or $\cos(\omega t)$, integrating over one period, and considering the steady state outputs that result when using $q$ different frequencies $\omega_1,\dots,\omega_q$, one can write the matrix equation
\begin{align} \label{forelate}
\frac{1}{\epsilon}\Gamma_1 &=     \pi  \Xi_1   \Upsilon_1 + O(\epsilon),
\end{align}
where
\begin{align*}
\Gamma_1 &=  \begin{bmatrix}
\omega_1 \int_0^{2\pi/\omega_1} y_{ss}(w_1,t) \sin(\omega_1 t) dt   \\ \omega_1 \int_0^{2\pi/\omega_1} y_{ss}(w_1,t) \cos(\omega_1 t) dt  \\ \vdots \\
\omega_q \int_0^{2\pi/\omega_q} y_{ss}(w_q,t) \sin(\omega_q t) dt   \\ \omega_q \int_0^{2\pi/\omega_q} y_{ss}(w_q,t) \cos(\omega_q t) dt
\end{bmatrix} \in \mathbb{C}^{2q \times 1},
\quad \Xi_1 = 
\begin{bmatrix}
s_{1,1}^0(\omega_1)   &   \dots &s_{1,M}^0(\omega_1) \\
c_{1,1}^0(\omega_1)   &   \dots &c_{1,M}^0(\omega_1) \\
&  \vdots  \\
s_{1,1}^0(\omega_q)   &   \dots &s_{1,M}^0(\omega_q) \\
c_{1,1}^0(\omega_q)   &   \dots &c_{1,M}^0(\omega_q) \\
\end{bmatrix} \in \mathbb{C}^{2q \times M},
\end{align*}
and $\Upsilon_1 = \begin{bmatrix} g_y^1 I_{1,A}^0 & \dots  &  g_y^M I_{M,A}^0  \end{bmatrix}^T \in \mathbb{C}^{M \times 1}$.  The individual terms of $\Upsilon_1$ can then be estimated, for instance by taking 
\begin{equation} \label{foestimate}
\Upsilon_1 = \frac{1}{\epsilon \pi} {\Xi_1}^\dagger \Gamma_1.
\end{equation}
A similar strategy was used in \cite{wils19cdc} to fit a linear model to the isostable reduced dynamical model.  As explained in \cite{wils19cdc}, isostable reduced models of the form \eqref{reducedmodinput} and \eqref{reducedoutput} can always be scaled so that $I_{n,A}^0$ can be chosen arbitrarily.  For this reason, it will be assumed that $I_{n,A}^0  = 1$ for all $n$ so that the relationship \eqref{forelate} can be used to identify the unknown $g^1_y, \dots, g^M_y$ terms.  


\subsection{Inferring Second Order Accurate Models for Single Output Systems} \label{secondsec}
Once the $O(\epsilon)$ dynamics have been determined, the steady state behavior of the $O(\epsilon^2)$ terms of \eqref{psiexp} can be used to determine the second order accurate terms of \eqref{taylgobs} and \eqref{tinput}.  To this end, recalling that $u(t) = \epsilon \sin(\omega t)$ one can substitute the asymptotic expansions \eqref{tinput} and \eqref{psiexp} into \eqref{reducedmodinput} and collect $O(\epsilon^2)$ terms to yield
\begin{equation} \label{firstord2}
\dot{\psi}_n^{(2)} = \lambda_n \psi_n^{(2)} + \sum_{k = 1}^M \psi_k^{(1)}  I_{n,A}^k \sin(\omega t).
\end{equation}
The steady state dynamics of each $\psi_k^{(1)}$ are given by \eqref{ord1sol} so that in steady state, \eqref{firstord2} can be written as
\begin{equation} \label{secondord2}
\dot{\psi}_{n,ss}^{(2)} = \lambda_n \psi_{n,ss}^{(2)} + \sum_{k = 1}^M    \big[   s_{1,k}^0(\omega) I_{k,A}^0 \sin(\omega t)  + c_{1,k}^0(\omega) I_{k,A}^0 \cos(\omega t)  \big]              I_{n,A}^k \sin(\omega t).
\end{equation}
Using the trigonometric product-to-sum identities from \eqref{prodid} in Appendix \ref{apxb}, one can rewrite \eqref{secondord2} as
\begin{align} \label{thirdord3}
\dot{\psi}_{n,ss}^{(2)} &= \lambda_n \psi_{n,ss}^{(2)} + \frac{1}{2} \sum_{k = 1}^M  \bigg[ I_{n,A}^k   I_{k,A}^0   \big( \sin(2 \omega t) c_{1,k}^0(\omega)  - \cos(2 \omega t) s_{1,k}^0( \omega) + s_{1,k}^0(\omega) \big) \bigg].
\end{align}
Because  \eqref{thirdord3} is a periodically forced linear equation, one can verify that the steady state solution is
\begin{align} \label{ord2sol}
\psi^{(2)}_{n,ss} &=   \frac{-1}{2(\lambda_n^2 + 4 \omega^2)}  \sum_{k = 1}^M \bigg[  \sin(2 \omega t)   I_{n,A}^k I_{k,A}^0 \big( \lambda_n  c_{1,k}^0(\omega)  +   2 \omega  s_{1,k}^0(\omega)  \big)   \nonumber \\
 &+ \cos(2 \omega t)   I_{n,A}^k  I_{k,A}^0 \big(  -\lambda_n s_{1,k}^0(\omega) + 2 \omega  c_{1,k}^0(\omega)  \big)  \bigg]   -\sum_{k = 1}^M \bigg[ \frac{s_{1,k}^0(\omega)   I_{n,A}^k   I_{k,A}^0 }{2 \lambda_n} \bigg] \nonumber \\
 &= \sum_{k = 1}^M \bigg[  I_{n,A}^k s_{2,n}^k(\omega) \sin(2 \omega t)  + I_{n,A}^k   c_{2,n}^k(\omega) \cos(2 \omega t) - \frac{s_{1,k}^0(\omega)   I_{n,A}^k   I_{k,A}^0 }{2 \lambda_n}   \bigg ],
\end{align}
where $s_{2,n}^k(\omega) \in \mathbb{C}$ and $s_{2,n}^k(\omega) \in  \mathbb{C}$ are defined appropriately.  Once again, considering the output given by \eqref{reducedoutput}, to leading order the $O(\epsilon^2)$ terms of the output equations are 
\begin{align} \label{fourthord2}
y_{ss}(\omega,t)-y_0 &=  O(\epsilon) +  \sum_{j = 1}^M \big[  \psi_j^{(2)} g_y^j   \big] + \sum_{j = 1}^M \sum_{k = 1}^j \big[ \psi_j^{(1)} \psi_k^{(1)} g_y^{jk}  \big]  + O(\epsilon^2) \nonumber \\
&= O(\epsilon) + \epsilon^2 \sum_{j = 1}^M  g_y^j \sum_{k = 1}^M \bigg[  I_{j,A}^k s_{2,j}^k(\omega) \sin(2 \omega t)  + I_{j,A}^k   c_{2,j}^k(\omega) \cos(2 \omega t) - \frac{s_{1,k}^0(\omega)   I_{j,A}^k   I_{k,A}^0 }{2 \lambda_j}   \bigg ]   \nonumber \\
& + \sum_{j = 1}^M \sum_{k = 1}^j \bigg[ g_y^{jk}\big(  s_{1,j}^0(\omega)  I_{j,A}^0   \sin(\omega t) +   c_{1,j}^0(\omega)  I_{j,A}^0  \cos(\omega t)  \big)  \nonumber \\
&\times \big(   s_{1,k}^0(\omega)  I_{k,A}^0   \sin(\omega t) +   c_{1,k}^0(\omega)  I_{k,A}^0  \cos(\omega t) \big)  \bigg] + O(\epsilon^3),
\end{align}
where the $O(\epsilon)$ terms are given by \eqref{xss}.  Once again, trigonometric product-to-sum identities from \eqref{prodid} can be used to rewrite  \eqref{fourthord2} according to
\begin{align} \label{fifthord2}
y_{ss}(\omega,t)-y_0 &= O(\epsilon) + \epsilon^2 \sum_{j = 1}^M  g_y^j \sum_{k = 1}^M \bigg[  I_{j,A}^k s_{2,j}^k(\omega) \sin(2 \omega t)  + I_{j,A}^k   c_{2,j}^k(\omega) \cos(2 \omega t) - \frac{s_{1,k}^0(\omega)   I_{j,A}^k   I_{k,A}^0 }{2 \lambda_j}   \bigg ]   \nonumber \\
& + \frac{1}{2} \sum_{j = 1}^M \sum_{k = 1}^j g_y^{jk}  \bigg[    - s_{1,j}^0(\omega)  I_{j,A}^0  s_{1,k}^0(\omega)  I_{k,A}^0 \cos(2 \omega t) +  c_{1,j}^0(\omega)  I_{j,A}^0  c_{1,k}^0(\omega)  I_{k,A}^0 \cos(2 \omega t) \nonumber \\
&+    s_{1,j}^0(\omega)  I_{j,A}^0  c_{1,k}^0(\omega)  I_{k,A}^0 \sin(2 \omega t)   +   c_{1,j}^0(\omega)  I_{j,A}^0   s_{1,k}^0(\omega)  I_{k,A}^0  \sin(2 \omega t)  \nonumber \\
&+    s_{1,j}^0(\omega)  I_{j,A}^0   s_{1,k}^0(\omega)  I_{k,A}^0  + c_{1,j}^0(\omega)  I_{j,A}^0   c_{1,k}^0(\omega)  I_{k,A}^0    \bigg] + O(\epsilon^3),
\end{align}
Notice that the $O(\epsilon^2)$ terms are a linear combination of sines, cosines, and constants.  Also recall that terms of the form $g_y^j$ and $I_{j,A}^0$ can be estimated independently for all $j$ using \eqref{forelate}.  Thus, the relationships \eqref{fifthord2} can be rewritten as
\begin{align} \label{sixthord2}
y_{ss}(\omega,t)-y_0 &= O(\epsilon) +   \epsilon^2 \sum_{j = 1}^M \sum_{k = 1}^M   I_{j,A}^k  \big( \beta_{j,k,1}(\omega)  \sin(2 \omega t) + \beta_{j,k,2} (\omega) \cos(2 \omega t) + \beta_{j,k,3}(\omega) \big) \nonumber \\
&+ \sum_{j = 1}^M  \sum_{k = 1}^j  g_y^{jk}  \big(   \beta_{j,k,4}(\omega)  \sin(2 \omega t) + \beta_{j,k,5}(\omega) \cos(2 \omega t) + \beta_{j,k,6}(\omega) \big) + O(\epsilon^3),
\end{align}
where $\beta_{j,k,1}(\omega), \dots,\beta_{j,k,6}(\omega) \in \mathbb{C}$ are defined so that \eqref{fifthord2} and \eqref{sixthord2} are identical.  Finally, recalling that the $O(\epsilon)$ terms from \eqref{sixthord2} are a linear combination of sines and cosines with frequency $\omega$, employing a similar approach that was used to identify \eqref{forelate} one can multiply both sides of \eqref{sixthord2} by either $\sin(2 \omega t)$, $\cos(2 \omega t)$, or $1$ and integrate over one period using $q$ different frequencies, $\omega_1,\dots,\omega_q$ to yield an equation of the form
\begin{align} \label{sorelate}
\frac{1}{\epsilon^2} \Gamma_2 &=   \pi  \Xi_2   \Upsilon_2 + O(\epsilon),
\end{align}
where
\begin{align*}
\Gamma_2 &=   \begin{bmatrix}
\omega_1 \int_0^{2\pi/\omega_1} y_{ss}(w_1,t) \sin(2 \omega_1 t) dt   \\ \omega_1 \int_0^{2\pi/\omega_1} y_{ss}(w_1,t) \cos(2 \omega_1 t) dt  \\   \omega_1 \int_0^{2\pi/\omega_1} y_{ss}(w_1,t)  dt  -x_0 \\  \vdots \\
\omega_q \int_0^{2\pi/\omega_q} y_{ss}(w_q,t) \sin(2 \omega_q t) dt   \\ \omega_q \int_0^{2\pi/\omega_q} y_{ss}(w_q,t) \cos(2 \omega_q t) dt  \\  \omega_q \int_0^{2\pi/\omega_q} y_{ss}(w_q,t)  dt  -x_0
\end{bmatrix} \in \mathbb{C}^{3q \times 1},
\end{align*}
$\Upsilon_2 \in \mathbb{C}^{ \left( \frac{3M^2+M}{2} \right) \times 1}$ is a vector containing all of the terms in the second order accurate expansion, i.e.,~ terms  of the form $I_{j,A}^k$  and $g_y^{jk}$, and $\Xi_2 \in \mathbb{C}^{3q \times  \left( \frac{3M^2+M}{2} \right) }$.  Much like for \eqref{forelate}, the second order accurate terms of the isostable reduced model \eqref{reducedmodinput} that comprise $\Upsilon_2$ can be estimated according to 
\begin{equation}\label{soestimate}
\Upsilon_2 = \frac{1}{\epsilon^2 \pi}  \Xi_2^\dagger \Gamma_2 .
\end{equation}

\subsection{Inferring Reduced  Models to Arbitrary Orders of Accuracy for Single Output Systems} \label{arbsec}
Sections \ref{firstsec} and \ref{secondsec} detail a strategy for inferring first and second order accurate reduced order isostable models from output data.  As illustrated in Appendix \ref{apxc}, it is possible to identify the terms of the expansions  \eqref{taylgobs} and \eqref{tinput} to arbitrary orders of accuracy using linear relationships of the form \eqref{arblin}.  Intuitively this is made possible by Lemma C.1 which states that under the application of an input $u(t) = \epsilon \sin(\omega t)$, the steady state output can be written as
\begin{equation} \label{restatelemma}
y_{ss}(\omega,t) - y_0 = \sum_{k = 0}^{j-1}  \bigg[  \tau_k(\omega) \sin(k \omega t) + \sigma_k(\omega) \cos(k \omega t)  \bigg] + \epsilon^j \bigg[ \tau_j(\omega) \sin(j \omega t) + \sigma_j(\omega) \cos(j \omega t) \bigg] + O(\epsilon^{j+1}),
\end{equation}
where $\tau_j(\omega) \in \mathbb{C}$ and $\sigma_j(\omega)  \in \mathbb{C}$.  In other words, the oscillatory components of the response with frequency $j \omega$ first appear in the $O(\epsilon^j)$ terms.   These specific terms can be isolated by multiplying \eqref{restatelemma} by either $\sin(j \omega t)$ or $\cos(j \omega t)$ and taking the integral over a single period.   Subsequently, relationships between the terms that comprise each $\tau_k(\omega)$ and $\sigma_j(\omega)$ can be derived to infer the unknown coefficients required for a  $j^{\rm th}$ order reduction as long as the coefficients of the first, second, $\dots$, $j-2^{\rm th}$ and $j-1^{\rm th}$ order reductions have been identified.  Explicit details of this fitting strategy are given in Appendix \ref{apxc}

As a final note, as the order of accuracy required increases, the complexity of the necessary terms grows rapidly.  For instance, the $O(\epsilon)$ terms of the steady state response are written succinctly in \eqref{xss}; the $O(\epsilon^2)$ terms of the steady state output \eqref{sixthord2} are much more complicated.  For reductions that are valid to higher orders of accuracy it advisable to compute the terms of each matrix $\Xi_j$, and $R_{j}$ from \eqref{arblin} using a symbolic computational package.

\subsection{Estimation of Linearized Eigenvalues}
Knowledge of the eigenvalues of the Jacobian from \eqref{locallin} (i.e.,~the decay rates of isostable coordinates) is necessary to identify the terms of reduced order models.  This can be accomplished in a two step procedure that involves first obtaining a coarse estimate of the eigenvalues from noisy data and subsequently refining these estimates using the relationship \eqref{forelate}.  To obtain a coarse estimation of the linearized eigenvalues as described below, it is assumed that low intensity background noise prevents the system from reaching the stable fixed point in the absence of stimulation.  If no background noise is present, one could alternatively apply a small magnitude, random input that simulates a white noise process.

\subsubsection{Coarse Estimation of Linearized Eigenvalues} \label{coarsesec}
A coarse approximation of the linearized eigenvalues can be obtained  using noisy, unperturbed model output.  This strategy is based on similar strategies implemented in \cite{wils20ddred} and \cite{wils20acc}.  This strategy represents a variation of the delay-embedding approach proposed in \cite{arba17}.  To begin, consider an isostable reduced  model \eqref{reducedmodinput} with output \eqref{reducedoutput}.  It will be assumed that the output $y(t)=O(\epsilon)$.  It will also be assumed that each eigenvalue is not repeated.  Consider a time series of output data taken every $\Delta t$ seconds.  This time series can be stored as a matrix $Y \in \mathbb{R}^{K,L}$ where
\begin{equation} \label{yform}
Y \equiv \begin{bmatrix} y_{t_0} & y_{t_1} & \dots & y_{t_{L-1}} \end{bmatrix},
\end{equation}
with $y_{t_i} \in \mathbb{R}^K$ defined according to  
\begin{equation}
y_{t_i} = \begin{bmatrix} y(K i \Delta t)-y_0 \\  y(K i \Delta t+ \Delta t)-y_0  \\ \vdots \\ y(K i \Delta t + \Delta t (K-1)) - y_0 \end{bmatrix}.
\end{equation}
As a preliminary step in the analysis, proper orthogonal decomposition (POD) \cite{holm96}, \cite{rowl17} can be employed to decompose the time series contained in $Y$ into a set of representative modes.  Assuming that $L > K$, the POD modes $\phi_i$, can be computed by finding the eigenvectors of the covariance matrix $Y Y^T$. The corresponding eigenvalues $\lambda_i^{\rm POD}$ give a sense of the importance of the temporal fluctuations captured by the associated POD mode.  In order to accurately represent the output data, a total of $N_{\rm POD}$  modes can be chosen so that $\sum_{j = 1}^{N_{\rm POD}} \lambda_j^{\rm POD}/\sum_{j = 1}^K \lambda_j^{\rm POD} \approx 1$.  Defining $\Phi \in \mathbb{R}^{K \times N_{POD}}$, such that $\Phi = \begin{bmatrix} \phi_1 & \dots &\phi_{N_{\rm POD}} \end{bmatrix} $, any column of $Y$ can be projected onto the new POD basis as
\begin{equation} \label{podout}
y_{t_i} - y_0 = \sum_{j = 1}^{N_{\rm POD}} \big( \phi_j p_j^i \big) + {\rm res}_i^{\rm POD} = \Phi \mu_i + {\rm res}_i^{\rm POD},
\end{equation}
where $\mu_i = \begin{bmatrix} p_j^i  & \dots & p_{N_{\rm POD}^i}  \end{bmatrix}^T$ is a vector of POD coefficients associated with $y_{t_i}$ and ${\rm res}_i^{\rm POD}$ is a small residual term that stems from the fact that some of the POD modes are truncated.  Simultaneously, considering \eqref{taylgobs}, to leading order $\epsilon$ the output is given by
\begin{equation} \label{isobasis}
y(t)-y_0 = \sum_{k = 1}^M \psi_k(t) g_y^k.
\end{equation}
Recalling that $\dot{\psi}_k = \lambda_k \psi_k$ so that $\psi_k(t) = \exp(\lambda_k t) \psi_k(0)$, one can write
\begin{equation} \label{isoout}
y_{t_i} - y_0 = P \Psi_i + {\rm res_i^{\rm ISO}}
\end{equation}
where $\Psi_i = \begin{bmatrix} \psi_1(K i \Delta t) & \dots & \psi_M(K i \Delta t) \end{bmatrix}^T$, $P \in \mathbb{C}^{K \times M}$ where the $k^{\rm th}$ row and $j^{\rm th}$ column are equal to $g_y^j \exp(\lambda_j \Delta t (k-1))$, and ${\rm res}_i^{ISO}$ is a small residual terms that results from truncating both the $O(\epsilon^2)$ terms and the most rapidly decaying isostable coordinates.  Assuming both residuals are negligible in \eqref{podout} and \eqref{isoout} and noticing that the left hand sides are identical, one can write
\begin{align}
\mu_i &= \Phi^T P \Psi_i,  \label{mueq}  \\
\Psi_i &= P^\dagger \Phi \mu_i, \label{psieq}
\end{align}
where $^\dagger$ is the Moore-Penrose pseudoinverse.  Above, equation \eqref{mueq} uses the fact that the POD modes are orthogonal and Equation \eqref{psieq} assumes that the columns of $P$ are linearly independent so that $P^\dagger P$ gives the identity matrix.  Letting $\Lambda = {\rm diag}(\exp (\lambda_1 K  \Delta t),   \dots , (\lambda_M K  \Delta t))$ note that $\Psi_i = \Lambda \Psi_{i+1}$  Using this information, along with \eqref{mueq} and \eqref{psieq}, one finds that 
\begin{equation}
\mu_{i+1} = \Phi^T P \Lambda P^\dagger \Phi \mu_i = A_\mu \mu_i.
\end{equation}
Additionally, as discussed in \cite{wils20ddred}, when residuals ${\rm res}_i^{\rm ISO}$ and ${\rm res}_i^{\rm ISO}$ are small enough that they are negligible, $A_\mu$ shares eigenvalues with $\Lambda$.  This allows each $\lambda_m$ term to be estimated from data by estimating the matrix $A_\mu$ according to 
\begin{equation}
A_\mu = X^+ {X^-}^\dagger
\end{equation}
where $X^+ = \begin{bmatrix} \mu_{1} & \dots  \mu_{{L-1}}  \end{bmatrix}$ and $X^- = \begin{bmatrix} \mu_{0} & \dots  \mu_{{L-2}}  \end{bmatrix}$.  Subsequently, each eigenvalue $\lambda_\mu$ of $A_\mu$ can be related to the decay rates of the isostable coordinates from \eqref{reducedmodinput} according to
\begin{equation}\label{evalues}
\lambda_j = \frac{\lambda_{\mu,j}}{K \Delta t}.
\end{equation}
Note that the accuracy  of \eqref{evalues} will be influenced by the magnitude of the residuals ${\rm res_i^{\rm ISO}}$ and ${\rm res_i^{\rm POD}}$ which  depend on factors such as noise, truncation of the POD modes, and truncation of the $O(\epsilon^2)$ terms from \eqref{isobasis}.  In practice, the above approach gives a coarse approximation of the eigenvalues that comprise the reduced order model \eqref{reducedmodinput}.

\subsubsection{Refinement of Eigenvalue Estimates} \label{refineestimate}
In some cases, the coarse eigenvalue estimates from the previous section can still yield accurate result when fitting reduced order models of the form \eqref{reducedmodinput} and \eqref{reducedoutput} using the methods described in Sections \ref{firstsec}-\ref{arbsec}.  In other cases, however, it can be useful to refine these estimates in conjunction with relationships derived in Section \ref{firstsec}.   To do so, recalling that each term $I_{n,A}^0$ is assumed to be equal to 1,  consider the relationship \eqref{forelate} that can be used to estimate the first order accurate terms of the reduction
\begin{equation} \label{hfn}
h(g_y^1, \dots, g_y^M, \lambda_1, \dots,\lambda_M) = \pi \Xi_1 \Upsilon_1 - \frac{1}{\epsilon} \Gamma_1.
\end{equation}
When $\lambda_1, \dots,\lambda_M$ are known exactly, it is possible to find values of each $g_y^1, \dots, g_y^M$ so that  $h(g_y^1, \dots, g_y^M, \lambda_1, \dots,\lambda_M) \approx 0$, for instance, by using Equation \eqref{foestimate}.   Supposing that only an estimate of $\lambda_1, \dots,\lambda_M$ is available (perhaps using the method from Section \ref{coarsesec}) with a corresponding estimate of  $g_y^1, \dots, g_y^M$,  denote  $Z = \begin{bmatrix} g_y^1 & \dots &  g_y^M & \lambda_1 & \dots \lambda_M \end{bmatrix}^T$ and let  $Z_0$ be the initial estimate.  One can update this estimate of $Z_0$ according to the Newton iteration $Z_1 = Z_0 + \Delta Z$ where 
\begin{equation} \label{newtonit}
\Delta Z =  - \bigg( \frac{\partial h}{\partial Z} \bigg|_{Z_0}  \bigg)^\dagger  h(Z_0).
\end{equation}
Above, $\Delta Z$ represents the least squares solution to $ \big( \frac{\partial h}{\partial Z} \big|_{Z_0}  \big) \Delta Z = -h(Z_0)$ and attempts to update the estimated parameters so that $ ||h(Z_0) +  \big( \frac{\partial h}{\partial Z} \big|_{Z_0}  \big) \Delta Z  ||   \approx  ||h(Z_0 + \Delta Z)|| \approx 0 $ where $||\cdot||$ denotes the Euclidean norm.   The update \eqref{newtonit} can be repeated until convergence is achieved.  Note that convergence is not guaranteed; like any Newton iteration, this strategy works best when starting with a good initial estimation of $Z$.  Additionally, if convergence is not achieved it may be necessary to change the number of isostable coordinates considered in the reduction.



\subsection{List of Steps Necessary for Obtaining Reduced Order Models}
The following summarizes a list of steps required to estimate the terms of \eqref{reducedmodinput} and \eqref{reducedoutput} with the expansions \eqref{taylgobs} and \eqref{tinput} used for $G_y$ and each $I_{k,A}$, respectively, computed to arbitrary order accuracy in the isostable coordinates.   
\begin{enumerate}[Step 1)]
\item Using the strategy given in Section \ref{coarsesec}, decide on the number of isostable coordinates to use in the reduction and obtain a coarse estimate for their unperturbed decay rates.  It is advisable to choose a number of POD modes, $N_{\rm POD}$ so that $\sum_{j = 1}^{N_{\rm POD}} \lambda_j^{\rm POD}/\sum_{j = 1}^K \lambda_j^{\rm POD} \approx 1$.   Choosing $N_{\rm POD}$ so that this ratio is closer to 1 will allow more of the dynamical behavior to be captured but will also increase the risk of overfitting.   Step 1 can be skipped if the eigenvalues are known {\it a priori}.  Step 1 can also be skipped by making an initial guesss for the number of isostable coordinates and their associated decay rates and updating these guesses with the Newton iteration from Step 4.  
\item Once the number of isostable coordinates has been determined, use a symbolic computational package to compute the terms of $\Xi_j$ for each $j$ up to the desired order of accuracy.  These matrices are defined as part of \eqref{forelate} and \eqref{sorelate} for first order and second order accuracy relationships, respectively.  Third order accuracy terms and higher are defined in \eqref{largeraccuracy}.  For third order accurate reductions and higher,  associated terms of the vectors $R_j$ as defined in \eqref{largeraccuracy} will also need to be computed symbolically.
\item For various frequencies $\omega$, apply an input $u(t) = \epsilon \sin(\omega t)$.  After transients die out, use the resulting steady state responses to determine the terms of each vector $\Gamma_j$ defined according to Equations  \eqref{forelate}, \eqref{sorelate}, and \eqref{gammaj1} to be used in determining the first, second, and higher order accuracy terms of the expansions, respectively.  As discussed in the examples from the following sections, it can be helpful to use larger values of $\epsilon$ for higher order accuracy terms, although, this is not absolutely necessary.
\item Apply the strategy from \eqref{refineestimate} to refine the eigenvalue estimates obtained in Step 1.   If the Newton iteration does not converge, one can repeat Steps 1 and 2 using a smaller number of isostable coordinates (or by using a different guess for the initial isostable coordinates).
\item Once a refined estimate for the eigenvalues associated with the decay rates of the isostable coordinates is obtained, this information can be used to evaluate $\Xi_1$ and obtain  an estimate for the first order accurate terms using \eqref{foestimate}.  This information can be used to compute a numerical approximation for $\Xi_2$ and subsequently estimate the second order accurate terms of the expansion according to \eqref{soestimate}.  In general, once the $j^{\rm th}$ order accuracy terms are computed, as discussed in Appendix \ref{apxc}, numerical approximations of $\Xi_{j+1}$ and $R_{j+1}$ can be computed and the relationship \eqref{arblin} can be used to estimate the $j+1^{\rm th}$ order terms.  This process can be repeated until the terms of \eqref{taylgobs} and \eqref{tinput} are computed to the desired orders of accuracy. 
\end{enumerate}

\subsection{Inferring Reduced Output Models for Multiple Outputs}  \label{multisec}
Systems with multiple outputs can be considered with straightforward modifications to the model identification algorithms.  To do so, suppose  that the state dynamics of a general differential equation of the form \eqref{maineq} can be characterized by an isostable coordinate framework according to \eqref{reducedmod}.  As in \eqref{reducedmodinput}, it will be assumed that a rank-1 input can be applied  of the form $U(t) =  A u(t)$, $A = \begin{bmatrix} A_1 & \dots & A_N  \end{bmatrix}^T$ and $u(t) \in \mathbb{R}$  yielding isostable dynamics of the form $\eqref{reducedmodinput}$.  Additionally, suppose that there are $N_y$ outputs of the form
\begin{align}
y_m &= y_{m,0} + G_{m,y}(\psi_1,\dots,\psi_M), \nonumber  \\
& m = 1,\dots,N_y.
\end{align}
Analogous to \eqref{taylgobs}, each $G_{m,y}$ can be Taylor expanded in terms of the isostable coordinates
\begin{equation}\label{taylgobsmulti}
G_{m,y}(\psi_1,\dots,\psi_M)  \approx \sum_{k = 1}^{M} \left[ \psi_k g_{m,y}^k \right] + \sum_{j = 1}^M \sum_{k = 1}^j \left[ \psi_j \psi_k g_{m,y}^{jk} \right]+ \sum_{i = 1}^M \sum_{j = 1}^i \sum_{k = 1}^j \left[ \psi_i \psi_j \psi_k g_{m,y}^{ijk}\right] + \dots,
\end{equation} 
for $m = 1,\dots,N_y$.  The isostable dynamics in the multiple output scenario can still be represented according to \eqref{tinput}; adding  extra outputs does not change the isostable coordinate dynamics.  Instead, each additional coordinate yields an additional $G_{m,y}(\psi_1,\dots,\psi_M)$ for which the terms of the expansion from \eqref{taylgobsmulti} must be fit.  

To identify a first order accurate model with multiple outputs, for example, each  $g_{m,y}^{j}$ can be fit to the output data.  To do so, the techniques from Sections \ref{firstsec}-\ref{arbsec} can be used to identify equations of the form \eqref{forelate}
\begin{align} \label{fomulti}
\frac{1}{\epsilon} \Gamma_{m,1} &= \pi \Xi_{1} \Upsilon_{m,1}, \nonumber \\
& m = 1,\dots,N_y,
\end{align}
where similar to the definition given as part of \eqref{forelate},  $\Gamma_{m,1}$ is a vector calculated considering the steady state behavior of  $y_m$ when applying input at various frequencies and $\Upsilon_1$ is a vector of terms characterizing the the first order accurate isostable dynamics from \eqref{tinput} and the first order accurate terms of the expansion of $G_{m,y}$ from  \eqref{taylgobsmulti}.  Note that the terms of $\Xi_1$ are identical for all output relationships and were defined as part of \eqref{forelate}.  As done in Section \ref{firstsec}, one can assume that $I_{n,A}^0 = 1$ for all $n$, so that the relationship \eqref{fomulti} can be used to identify $g_{m,y}^1, \dots, g_{m,y}^M$ for all $m$. 

Once the first order terms are identified, the second order accurate terms can be inferred by identifying equations of the form \eqref{sorelate} for each output
\begin{align} \label{somulti}
\frac{1}{\epsilon^2} \Gamma_{m,2} &= \pi \Xi_{2} \Upsilon_{m,2}, \nonumber \\
& m = 1,\dots,N_y,
\end{align}
 where $\Gamma_{m,2}$ and $\Upsilon_{m,2}$ are defined analogously to $\Gamma_{2}$ and $\Upsilon_{2}$ from \eqref{fomulti} for each output $y_m$.  Note that $\Xi_2$ was defined as part of \eqref{sorelate}.    All relationships of the form \eqref{somulti} for $m = 1,\dots,N_y$ can be used to infer the second order terms of each ${I_{n,A}}$ from \eqref{tinput}  and of each $G_{m,y}$ from \eqref{taylgobsmulti}. In general, once the $j^{\rm th}$ order accurate terms have been estimated, a similar strategy can be used to infer the $j+1^{\rm th}$ order terms of the expansion by identifying equations of the form \eqref{largeraccuracy} for each output
 \begin{align}
 \Gamma_{m,j+1} = \pi \epsilon^{j+1} \big( \Xi_{j+1} \Upsilon_{m,j+1} + R_{m,j+1}   \big) + O(\epsilon^{j+2}) ,
 \end{align}
  where $\Gamma_{m,n+1}$, $\Upsilon_{m,j+1}$, and $R_{m,j+1}$ are defined analogously to $\Gamma_{j+1}$, $\Upsilon_{j+1}$, and $R_{j+1}$ from \eqref{largeraccuracy} for each input $y_m$ and $\Xi_{j+1}$ is defined as part of Equation \eqref{largeraccuracy}.

As a final consideration, note that the number of terms to fit grows in proportion to the number of outputs.  In situations where  there are a large number of outputs, it can be useful to first apply a reduction strategy such as POD to the outputs and consider a subset of the most important POD coefficients to be the system output rather than fitting coefficients for the full set of outputs.  This is intuitively similar to Galerkin projection methods for obtaining numerical solutions of partial differential equations, for example, as performed in \cite{camp05}.  This suggested approach will be illustrated in the example from Section \ref{burgsec}.

\section{Examples} \label{exampsec}

\subsection{Reduced Order Modeling of a Simple Dynamical System with Additive Noise} \label{simplesec}
As a preliminary illustration of the proposed reduced order model identification strategy, consider the simple dynamical system from \eqref{exampfn} with white noise added to the $x_1$ variable
\begin{align} \label{exampfnnoise}
\dot{x}_1 & =  \mu x_1 + u_{x_1}(t) + \sqrt{2 D} \eta(t), \nonumber \\
\dot{x}_2 & = \lambda(-x_1 + x_2 + x_1^2 + x_1^3),
\end{align}
where $\lambda = -1$, $\mu = -0.05$, and $\eta(t)$ is an independent and identically distributed zero-mean white noise process with intensity $D = 0.0005$.  The model output is taken to be $y(t) = x_2(t)$.  Note that the only difference between \eqref{exampfn} and the above equation is the addition of noise.   In Section \ref{reducedexamp} a reduced order model was determined with knowledge of the underlying model Equations \eqref{exampfn}.  Here, a reduced order model of the form \eqref{reducedmodinput} with output of the form \eqref{reducedoutput} will be inferred strictly from measurements of model output without any assumed knowledge of the underlying system dynamics.  To do so, \eqref{exampfnnoise} is simulated for 20,000 time units with a time step of 0.1 using the stochastic simulation algorithm presented in \cite{hone92}.  The resulting data is stored in a matrix  of the form \eqref{yform} using $K = 100$ and POD is performed on the resulting data.  Using a single POD mode, $\lambda_1^{\rm POD}/ \sum_{j = 1}^K (\lambda_j^{\rm POD}) = 0.87$.  The coarse eigenvalue estimation strategy detailed in  Section \ref{coarsesec} results in an initial eigenvalue estimate of $\lambda_1 = -0.0322$.

\begin{figure}[htb]
\begin{center}
\includegraphics[height=4in]{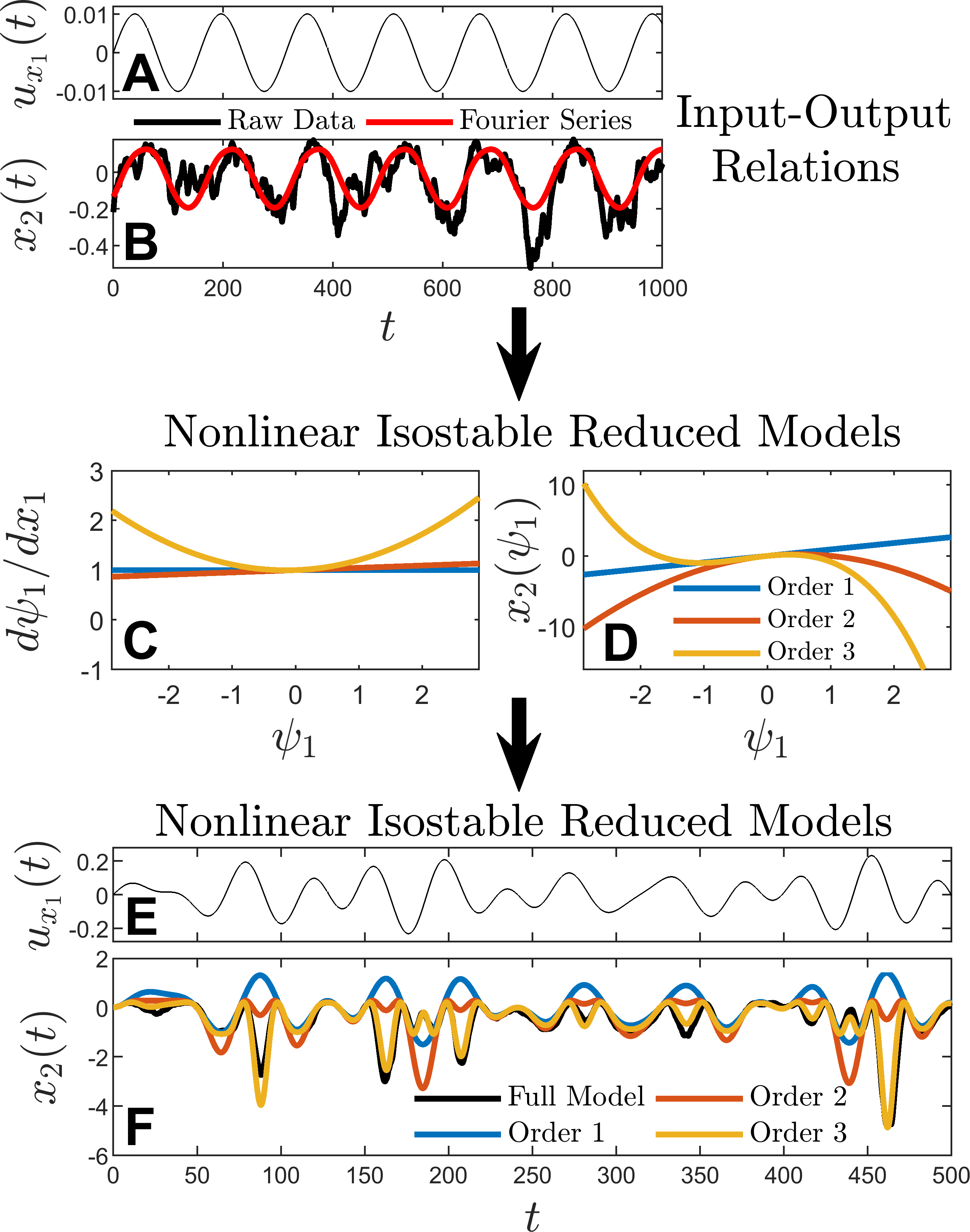}
\end{center}
\caption{Results from applying the model inference strategy described in Section \ref{fittingsection} to the simple model \eqref{exampfnnoise}. Panel A shows an example of a sinusoidal input with the black line in panel B giving the corresponding model output after allowing the transient dynamics to decay.  Fourier coefficients are computed for the first, second, and third harmonics by averaging the output over 100 cycles and the associated Fourier series representation is shown as a red line.  Using inputs of various frequency and magnitude, the corresponding outputs are measured and this information is used to identify reduced order models valid to various orders of accuracy as described in Section \ref{fittingsection}.  Curves from panels C and D characterize the response of the isostable coordinate to inputs and the output as a function of the isostable coordinate, respectively.  Notice that these curves are similar to those shown in panels C and D of Figure \ref{exampfig} which were obtained with full knowledge of the underlying dynamical equations.  The input shown in panel E is applied to the model equations \eqref{exampfnnoise} and the output in panel F is compared to the output predicted from the reduced order models valid to different orders of accuracy.   }
\label{inferredresults}
\end{figure}

Next, the stimulus $u_{x_1}(t) = \epsilon \sin(\omega t)$ is applied taking  $\omega \in \{0.02,0.025,0.03,0.035,0.04\}$.  Using $\epsilon = 0.01$, the steady state output is averaged over 100 cycles in order to compute the terms of $\Gamma_1$ according to \eqref{forelate}.  An example input-output relationship in steady state is shown in Panels A and B of Figure \ref{inferredresults}.  The terms of $\Gamma_2$ (from \eqref{sorelate}) and $\Gamma_3$  (from \eqref{gammaj1}) are computed taking $\epsilon = \sqrt{0.01}$ and $\sqrt[3]{0.01}$, respectively.  Choosing $\epsilon$ in this manner helps to emphasize the contributions of the higher order terms when estimating the associated coefficients of the expansions \eqref{taylgobs} and \eqref{tinput}. Auxiliary matrices  $\Xi_1$, $\Xi_2$, $\Xi_3$ and $R_3$ used to fit the first, second, and third order coefficients of the reduction are computed symbolically.  After obtaining a preliminary estimate of $g_y^1$, the Newton iteration suggested by \eqref{newtonit} is employed to refine the estimate of the principal eigenvalue of $\lambda_1 = -0.0462$ (which is relatively close to the true value of $\lambda_1 = -0.05$ for the noiseless system \eqref{exampfn}).  Second and third order terms are computed in succession using \eqref{soestimate} and \eqref{arblin}, respectively.   Panels C and D show estimates of $d\psi_1/d x_1$ and $x_2(\psi_1)$ obtained from first, second, and third order accurate estimates  of the reduced isostable dynamics for models of the form \eqref{reducedmodinput} with output \eqref{reducedoutput}.  These estimates are similar to results from panels C and D from Figure \ref{exampfig} that were obtained with knowledge of the full model equations.  The resulting estimates of the reduced order model parameters can be used to predict the response to other inputs.  Starting from an initial condition at $x_1 = x_2 = 0$, the model \eqref{exampfnnoise} is simulated taking $u_{x_1}(t) = 0.08(\sin(0.1t) + \sin(0.15 t) + \sin(0.17t))$ as shown in panel E.  The full model output is compared to the outputs from reduced models of various orders of accuracy.  For results shown here, the reduced models do not contain any noise terms.   Much like the results presented in Figure \ref{exampfig}, the third order accurate reduced model accurately reflects the output from the full simulations.   The first and second order accurate models perform decently, especially in moments when the magnitude of the input is small, but are not sufficient to replicate the full system behavior.     Note here that $u_{x_1}(t)$ is comprised of higher frequency sinusoids than those used to infer the reduced order models.  Additionally, the applied input is periodic, but has a period that is significantly larger than the sinusoids used to infer the reduced order model.  

\FloatBarrier

\subsection{Spike Rates of Neural Populations} \label{spikesec}
Next, the proposed reduced model inference strategy is applied to a more complicated model describing the spike rates of a population of synaptically coupled neurons from \cite{rubi04}.  The model equations are
\begin{align} \label{multimodel}
C \dot{V}_i &= -I_{\rm L}(V_i) - I_{\rm Na}(V_i,h_i) - I_{\rm K}(V_i,h_i) - I_{\rm T}(V_i,r_i)  \nonumber \\
& \quad + I^{b}_i - \frac{g_{\rm syn}}{N} \sum_{j =1}^N s_j (V_i-E_{\rm syn}) + \sqrt{2 D} \eta_i(t) + u(t), \nonumber \\
\dot{h}_i &= (h_\infty(V_i) - h_i)/\tau_h(V_i), \nonumber \\ 
\dot{r}_i  &= (r_\infty(V_i)-r_i)/\tau_r(V_i), \nonumber \\
\dot{s}_i &= \frac{a(1-s)}{1+\exp(-(V_i-V_T)/\sigma_T)} - b s_i.
\end{align}
Here, $N = 1000$ gives the total number of neurons in the population, $V_i$, $s_i$, $h_i$, and $r_i$ represent the transmembrane voltage, a synaptic variable, and two gating variables assigned to neuron $i$,  respectively, the conductance $g_{\rm syn}$ sets the coupling strength, $E_{\rm syn} = -100$ mV is the reversal potential of the neurotransmitter resulting in inhibitory coupling, $u(t)$ is an injected current common to all neurons, $C = 1 \mu {\rm F}/{\rm cm}^2$ is the membrane capacitance, $\sqrt{2D}\eta_i(t)$ is an independent and identically distributed zero-mean white noise process with intensity $D = 1$, and $I^b_i$ is the baseline current of neuron $i$ and is drawn from a normal distribution with a mean of 5 and a variance of 1 $\mu {\rm A}/{\rm cm}^2$.  The characteristics of the synaptic current are determined by the parameters $a = 3, V_t = -20 {\rm mV}, \sigma_T = 0.8 {\rm mV},$ and $ \beta = 1$.  The reader is referred to \cite{rubi04} for a full explanation of the remaining functions that determine the ionic currents $I_{\rm L}$, $I_{\rm Na}$, $I_{\rm K}$, $I_{\rm T}$ and the behavior of the gating variables.  The observable for \eqref{multimodel} is taken to be the firing rate $R(t)$ where neuron $i$ is defined to fire at the moment $V_i$ crosses -25 mV with a positive slope.   The firing rate defined over a sliding window
\begin{equation} \label{neuraloutput}
FR(t)  = \frac { \text{Total Number of Neurons That Fired in the Interval [t-W,t]} }{W},
\end{equation}
where $W = 1.5$ ms is the width of the window.

\begin{figure}[htb]
\begin{center}
\includegraphics[height=2.2in]{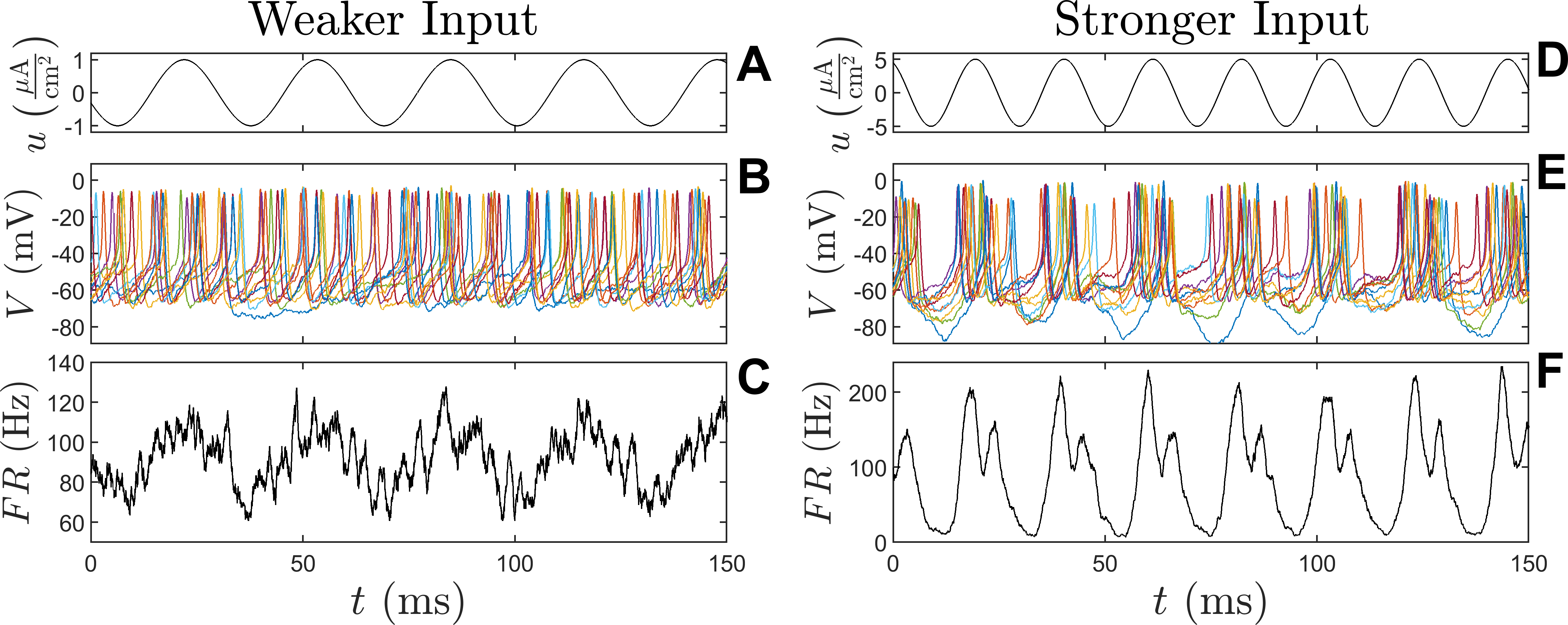}
\end{center}
\caption{ Panel A shows a sinusoidal inputs to \eqref{multimodel} with the resulting voltage traces from a subset of neurons and firing rates shown in panels B and C, respectively.  The same information is shown in panels D-F using a larger magnitude input with a different frequency.  For small inputs, $FR(t)$ appears somewhat sinusoidal.  When larger magnitude inputs are used the resulting firing rates are non-sinusoidal indicating the importance of nonlinear models to characterize the output in response to inputs.}
\label{showresults}
\end{figure}

Figure \ref{showresults} shows an example of the typical behavior of the model \eqref{multimodel} with output \eqref{neuraloutput} in response to input.  Panel A (resp.,~D) show an example of a weak (resp.~strong) sinusoidal input.  Voltage traces from 10 representative neurons within the population are shown in panel B (resp.~E) with the corresponding firing rate shown in panel C (resp.~F).  Towards obtaining a reduced order model to characterize the input-output relationships of this model, \eqref{multimodel} is simulated taking $\Delta t = 0.015$ ms for a total of 5000 ms.  The resulting data is stored in a matrix of the form \eqref{yform} taking $K = 50$ and POD is performed on the resulting data.  Using two POD modes yields  $ \sum_{j = 1}^2 ( \lambda_j^{\rm POD})/ \sum_{j = 1}^K (\lambda_j^{\rm POD}) = 0.93$.  Using the coarse eigenvalue estimation strategy from Section  \ref{coarsesec} in conjunction with these two POD modes  yields two complex conjugate eigenvalues $\lambda_{1,2} = -.65 \pm 1.65 i$.

A reduced order model is obtained by considering the firing rate in response to sinusoidal input $u(t) = \epsilon \sin(\omega t)$ with $\omega \in \{0.3,0.4,0.5,0.6,0.7, 0.8,0.9 \}$.  Taking $\epsilon  = 0.1 \mu {\rm A}/\mu {\rm F}$, initial transients are allowed to decay and the steady state output is averaged over 700 ms to compute the terms of $\Gamma_1$ from \eqref{forelate}.   In a similar manner, $\epsilon  = 1 \mu {\rm A}/\mu {\rm F}$ is used to compute the terms of  $\Gamma_2$ from \eqref{sorelate}.  Terms of the auxiliary matrices $\Xi_1$ and $\Xi_2$ from \eqref{forelate} and \eqref{sorelate}, respectively, are computed symbolically.  After obtaining preliminary estimates of the first order terms, the Newton iteration from \eqref{newtonit} is employed to obtain a refined estimate of $\lambda_{1,2} = -.19 \pm 0.68 i$ for the principal eigenvalues.  Second order terms are subsequently computed using \eqref{soestimate}.   Third order accurate terms were also considered, however, the resulting reduced order models did not perform better than the second order accurate models, and the results are not shown.

\begin{figure}[htb]
\begin{center}
\includegraphics[height=3.4 in]{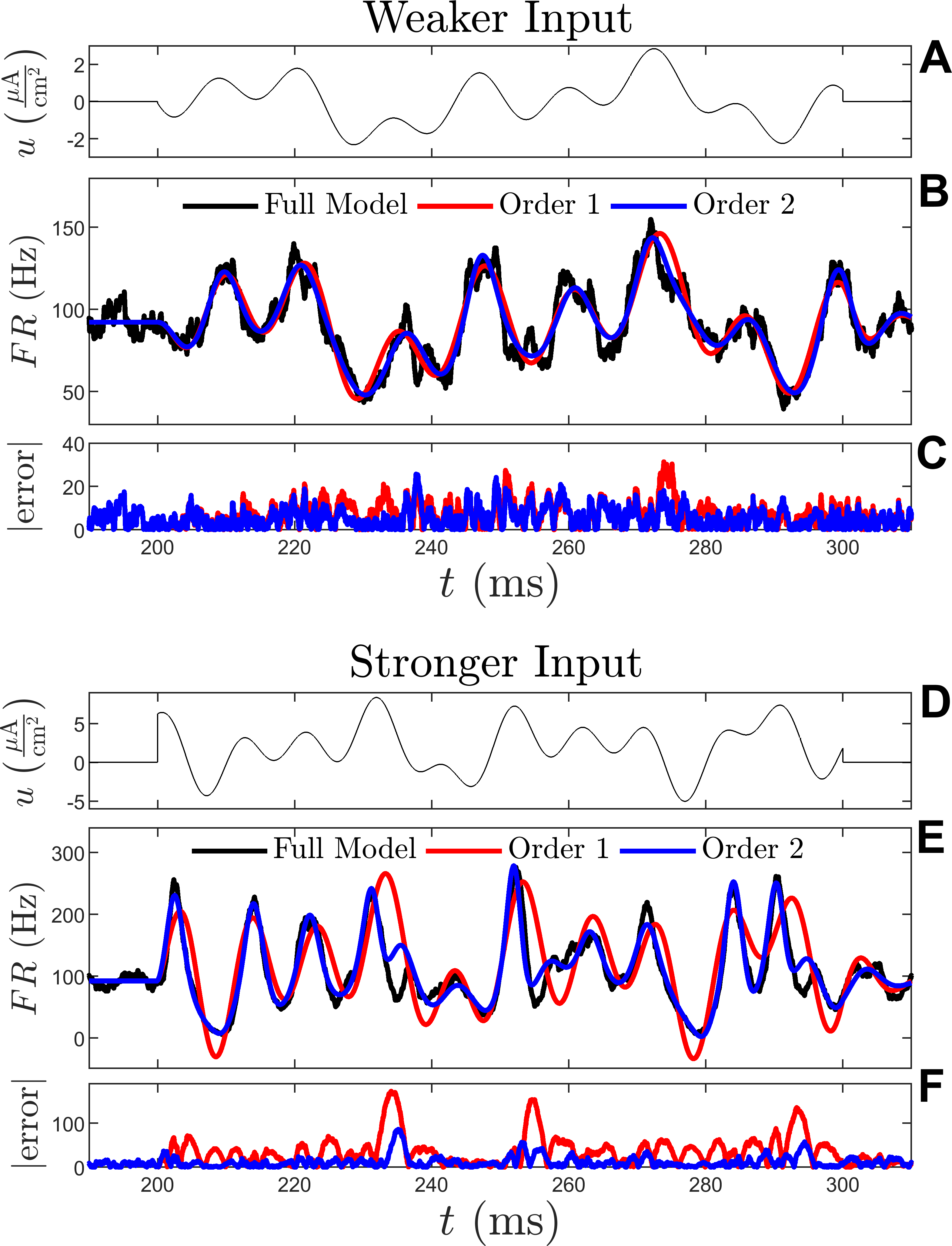}
\end{center}
\caption{ For the relatively weak input shown in panel A, the output of the first and second order accurate reduced order models in panel B (red and blue lines, respectively) closely resembles the output from the full model equations (black line).  Panel C shows the corresponding absolute error between the full and reduced order simulations.  Analogous results are shown for a larger magnitude input in panels D-F.  While both the first and second order accurate models perform well, the second order accurate model performs significantly better when larger magnitude inputs are considered.}
\label{neuroninput}
\end{figure}

The predictive capability of the resulting reduced order model is considered with two different inputs shown in Figure \ref{neuroninput}. In both the full \eqref{multimodel} and reduced model equations, no input is applied for the first 200 ms of simulation to allow for any transient behavior to die out.  Subsequently, a small magnitude input $u(t) = \sin(0.1 t) + \cos(0.23 t) + \sin(-.49 t)$ is applied for 100 ms and is shown in panel A.  Resulting outputs from the full and reduced models are shown in panel B with the absolute error between the reduced and full models shown in panel C.   Analogous results are shown  in panels D-F for a larger magnitude input $u(t) = 2.5(\sin(0.2 t) + \cos(0.35 t) + \sin(0.63t))$.  Note that each neuron from the full model equations \eqref{multimodel} has a white noise term added to the voltage equations while the reduced order models are deterministic -- this contributes to a jagged appearance of the full model outputs.  First and second order reduced models have comparable performance when weaker input is applied.  As the inputs become stronger, the performance of the first order accurate model begins to degrade while the second order accurate model still performs well.  Note that while the terms of the reduced order models were determined by considering the steady state behavior in response to purely sinusoidal forcing, they accurately replicate the transient response to the inputs considered in panels A and D (which are not periodic on the time scales considered).

\FloatBarrier

\subsection{One-Dimensional Burgers' Equation} \label{burgsec}
As a final example, consider the 1-D Burgers' equation
\begin{equation} \label{burgmod}
\frac{\partial w}{\partial t} = \frac{1}{{\rm Re}} \frac{\partial^2 w}{\partial x^2} - w \frac{\partial w}{\partial x},
\end{equation}
Here, $w$ is the state on the spatial domain $x \in [0,1]$, and ${\rm Re}$ is a viscosity term that is analogous to the Reynolds number from the Navier-Stokes equation.  The boundary condition at $x = 0$ is defined to be $w_L(t)= 0.3 + u(t)$ where 0.3 is the nominal value and $u(t)$ serves as a time-varying input.   The boundary condition at $x = 1$ is held fixed at a value of 0.3.    Because of its similar structure to the Navier-Stokes equation, the Burgers' equation is often  often used as a test bed for model reduction techniques in fluid flow applications  \cite{page18}, \cite{kutz18}.  

\begin{figure}[htb]
\begin{center}
\includegraphics[height=2.3in]{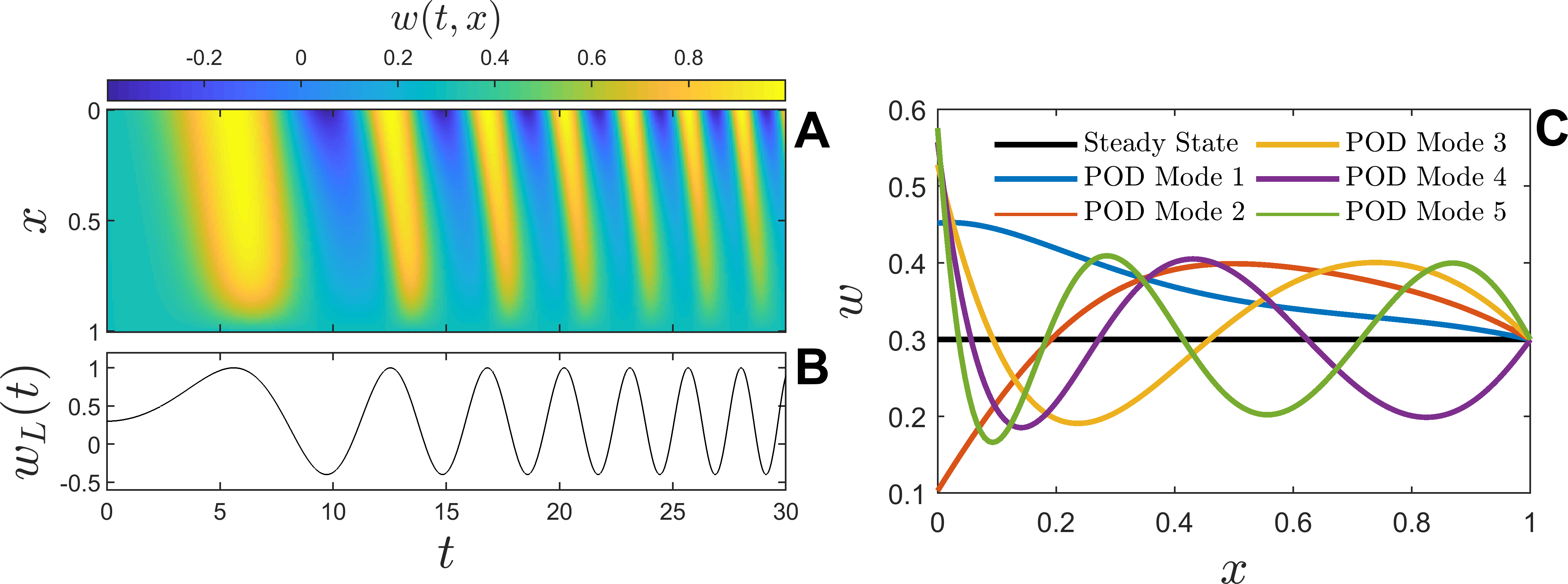}
\end{center}
\caption{Panel A shows the response of the 1-D Burgers' equation to the time-varying left boundary condition shown in panel B.  POD modes are extracted from the model output and plotted in Panel C as colored lines.  The steady state solution when $u(t) = 0$ is shown for reference as a black line. }
\label{burgmodes}
\end{figure}

In this example, ${\rm RE} = 10$ is used and the model output is taken to be the state on the entire spatial domain.  Towards determining appropriate POD modes to represent the model outputs, an input $w_L(t) = 0.3 + 0.7 \sin(t^2/20)$ is applied  for $t = 100$ time units in order to excite system modes over wide range of frequencies.   The model response to the time-varying output are shown in panels A and B of Figure \ref{burgmodes}, respectively.   POD is performed on the resulting snapshots taken at $\Delta t = 0.01$ increments with the spatial domain discretized into 152 equally spaced elements.  Taking 5 POD modes yields $ \sum_{j = 1}^5 ( \lambda_j^{\rm POD})/ \sum_{j = 1}^{152} (\lambda_j^{\rm POD}) = 0.99994$ and captures the majority of the deviations from the steady state solutions.  These resulting POD modes are shown in Panel C and are ranked in descending order of their importance.

A reduced order model using three isostable coordinates will be used with an initial guess of $\lambda_1 = -1$, $\lambda_2 = -2$, and $\lambda_3 = -3$.  As discussed in Section \ref{multisec}, the POD coefficients associated with the five POD modes from Figure \ref{burgmodes} are taken as the system output.   A reduced order model is obtained by considering the steady state model output (i.e.,~the time course of the POD coefficients in steady state) in response to input $u(t) = \epsilon \sin(\omega t)$ with $\omega = \{0.1,0.2,\dots,2.0\}$.  Taking $\epsilon = 0.05$ (resp.~0.5), the terms of $\Gamma_1$ (resp.~$\Gamma_2$)  from \eqref{forelate} (resp.~\eqref{sorelate}) are computed.  Terms of auxiliary matrices $\Xi_1$ and $\Xi_2$ from \eqref{forelate} and \eqref{sorelate}, respectively, are computed using a symbolic computational package.

The Newton iteration from \eqref{newtonit} is applied to obtain a refined estimate of the principal eigenvalues.  This Newton iteration is performed as described in Section \ref{refineestimate} using the first POD mode as the output.    This iteration converges to $\lambda_1 = -1.22$, $\lambda_2 = -4.62$, and $\lambda_3 = -32.06$.  The eigenvalues that result from the Newton iteration are robust to changes in the initial guess for the eigenvalues.   First and second order terms of the reduced functions (i.e.,~the functions that characterize the isostable dynamics and the isostable-to-output relationships for each POD coefficient) are then computed using relationships of the form \eqref{fomulti} and \eqref{somulti}. Reduced order models beyond second order accuracy are not considered for this example.   


\begin{figure}[htb]
\begin{center}
\includegraphics[height=3in]{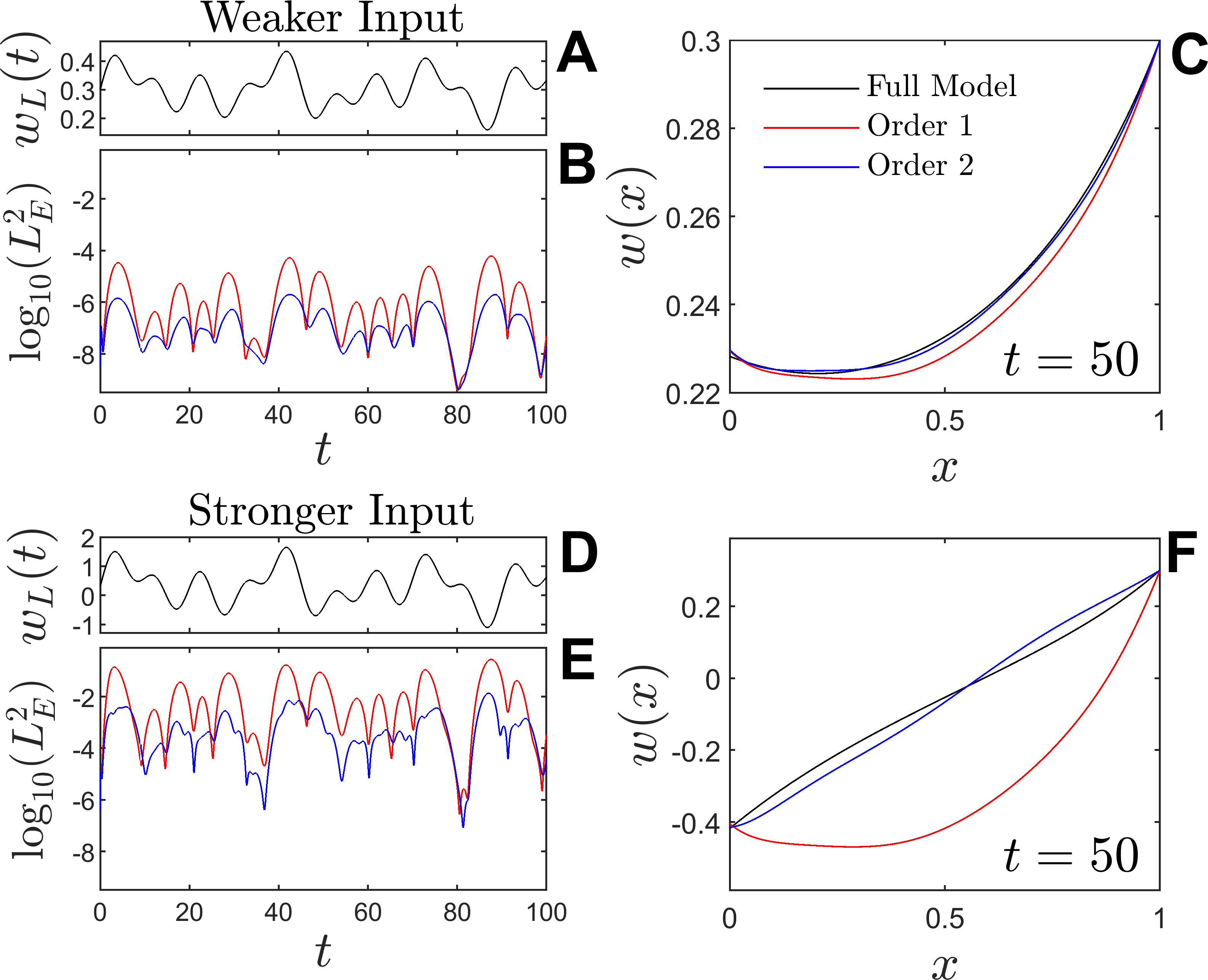}
\end{center}
\caption{The input $w_L(t)$ shown in panel A is applied to the full model eqref{burgmod} as well as the first and second order reduced models.  Panel B shows the resulting $L^2$ error between full and reduced order models.  Panel C gives a snapshot of the solution profiles for each model at $t = 50$.  Analogous results are shown in panels D-F when the stronger input from panel D is applied.  The second order model output has about two orders of magnitude less error than the first order model.   These differences are more pronounced when larger magnitude inputs are used as evidenced by panel F. }
\label{inputresult}
\end{figure}

Figure \ref{inputresult} compares the outputs predicted by the resulting reduced order models with isostable dynamics of the form \eqref{reducedmodinput} and model outputs \eqref{taylgobsmulti} to the full model dynamics from \eqref{burgmod}.    Here, $w_{\rm red}(x,t) = \sum_{j = 1}^5 \phi_j(x) p_j(t)$ where  $\phi_j(x)$ and $p_j(t)$ are the $j^{\rm th}$ POD modes and POD coefficients, respectively.   Starting from steady state the left boundary condition is taken to be  $w_L(t) = 0.3 + 0.05(\sin(0.2 t) + \sin(0.35 t) + \sin(0.63 t))$.  The logarithm of the $L^2$ error defined as $L^2_E(t) = \int_0^1(w_{\rm red}(x,t)-w_{\rm full}(x,t))^2 dx$ is shown in panel B, where $w_{\rm red}$ and  $w_{\rm full}$ comes from the reduced and full model simulations, respectively. The solution profiles at $t = 50$ are shown in panel C.  Panels D-F give analogous results using a larger magnitude input $w_L(t) = 0.3 + 0.5(\sin(0.2 t) + \sin(0.35 t) + \sin(0.63 t))$.  In both cases, the second order accurate models result in $L^2$ errors that are approximately two orders of magnitude smaller than the $L^2$ errors from the first order accurate models.  The differences in errors between these reduced order models are particularly pronounced  when using stronger inputs.

\FloatBarrier

\section{Conclusion} \label{concsec}
In applications where the underlying model equations are known, it is relatively straightforward to to numerically compute isostable-based reduced order models directly from the dynamical equations.  However, in most experimental applications, model equations are not known precisely and data-driven techniques must be used instead.      In this work, a general, data-driven strategy is developed for inferring nonlinear isostable reduced models of the form \eqref{reducedmodinput} and \eqref{reducedoutput}.

The strategy proposed in this work leverages an isostable-based coordinate system that characterizes the behavior of the level sets of the slowest decaying Koopman eigenfunctions.    Using this isostable coordinate basis as a foundation, one can show that in response to inputs of the form $u(t) = \epsilon \sin(\omega t)$, steady state model outputs can be written in terms of a Fourier series expansion with a fundamental frequency $\omega$.  Furthermore, components of the output with frequency $j \omega$ are shown to be $\mathcal{O}(\epsilon^j)$ terms.  These key insights are used to develop the proposed model inference strategy discussed in Sections \ref{firstsec}-\ref{arbsec}.  Using this strategy, it is possible to consider the steady-state output in response to sinusoidal input at various frequencies to infer high-accuracy reduced order models of the form \eqref{reducedmodinput} and \eqref{reducedoutput}.   The proposed model inference strategy is illustrated for three example systems.   In each example considered, a linear model that only considers the first order terms in the expansions of \eqref{taylgobs} and \eqref{tinput} performs reasonably well when small magnitude inputs are considered.  However, when considering larger magnitude inputs, nonlinear reduced order models are necessary to  capture the dynamical behavior.  



 While it is technically possible to infer the terms of the expansions \eqref{taylgobs} and \eqref{tinput} to arbitrary orders of accuracy in the isostable coordinates using the proposed method, the robustness of the model fitting procedure suffers as the desired order of accuracy increases.  Indeed, considering the structure of \eqref{arblin}, the output measurements that comprise $\Gamma_j$ are multiplied by $1/\epsilon^j$; this amplifies the influence of any errors or uncertainties associated with the measurements of the steady state output.  This effect can be partially mitigated by increasing the amplitude of the sinusoidal inputs as higher order accuracy terms are considered (this strategy was used in the examples considered in this work).  As an additional consideration, the  $j^{\rm th}$ order accuracy terms depend on the $j-1^{\rm th}$ order accuracy terms, so any errors in the computation of lower order accuracy terms will be propagated when computing higher order accuracy terms.  For the simple model \eqref{exampfnnoise} from Section \ref{simplesec} the third order terms of \eqref{taylgobs} and \eqref{tinput} were inferred reliably.  In the more complicated neural spike rate example From Section \ref{spikesec}, the performance of the reduced models was significantly improved by the addition of the second order terms, but  the third order terms did not improve the performance -- in this example it is possible that the contribution from the third order terms is simply negligible.  It would be worthwhile to investigate other statistical techniques to infer the reduced order terms from the relationships \eqref{forelate}, \eqref{sorelate}, and \eqref{largeraccuracy} for the first, second, and arbitrary order accuracy terms, respectively, perhaps involving either approximate Bayesian computation \cite{toni09}, \cite{sunn13} or deep learning approaches \cite{lusc18}, \cite{yeun19}.
 

Future work will investigate the  possibility of extending the proposed methods for use in oscillatory dynamical systems.  The `direct method' \cite{neto12}, \cite{izhi07} is a well-established technique that can be used to determine reduced order equations of the form $\dot{\theta} = \omega + Z(\theta)u(t)$ from experimental data, where $\theta \in \mathbb{S}^1$ is the phase, $Z(\theta)$ captures the effect of an input $u(t)$ on the phase, and $\omega$ is the nominal frequency of oscillation.  The reduced order dynamical models considered in this work from Equation \eqref{reducedmod} are similar in structure to the models that were investigated in \cite{wils20highacc} for oscillatory dynamical systems.  It would be of interest to investigate whether the data-driven model inference techniques proposed here can be adapted for use in applications that involve limit cycle oscillators.  

\section*{Acknowledgment}
This material is based upon work supported by the National Science Foundation Grant CMMI-2024111.

\section*{Data Availability Statement}
The data that support the findings of this study are available from the corresponding author upon reasonable request.

\begin{appendices}

\section{Asymptotic Expansions of Isostable Reduced Models} \label{apxa}

\renewcommand{\thetable}{A\arabic{table}}  
\renewcommand{\thefigure}{A\arabic{figure}} 
\renewcommand{\theequation}{A\arabic{equation}} 
\setcounter{equation}{0}
\setcounter{figure}{0}

In \cite{wils20highacc}, a strategy was devised for computing asymptotic expansions of the state and the gradient of the isostable coordinates for periodic orbit attractors in terms of the isostable coordinates themselves.    A nearly identical approach can be used to compute the individual terms of \eqref{taylG} and \eqref{taylI} for stable fixed point attractors.  The following results can be obtained by following the same arguments as those presented in \cite{wils20highacc}, but taking the attractor to be a fixed point  instead of a periodic orbit.  


First, assuming $U(t) = 0$,  the individual terms of $G(\psi_1,\dots,\psi_M)$ from \eqref{taylG} can be found by taking the time derivatives of $\Delta x = G(\psi_1,\dots,\psi_M)$ yielding
\begin{align} \label{timeexpand}
\frac{d \Delta x}{dt} &=  \sum_{k = 1}^{M} \left[   g^k \lambda_k \psi_k  \right] + \sum_{j = 1}^M \sum_{k = 1}^j \left[  g^{jk} (\lambda_j+\lambda_k) \psi_j \psi_k  \right] +\dots
\end{align}
Above, the relationship $\dot{\psi}_k = \lambda_k \psi_k$ is used when taking time derivatives using the product rule.  Letting $F(x) = \begin{bmatrix} f_1 & \dots & f_N \end{bmatrix}^T$, a direct expansion of \eqref{maineq} in terms of $\Delta x$ yields
\begin{align} \label{taylfn}
\frac{d \Delta x}{dt} &= J   \Delta x  +  \begin{bmatrix}  \sum_{i = 2}^\infty \frac{1}{i !} \left[ \overset{i}{\otimes}  \Delta x^T   \right]{\rm vec}( f_1^{(i)})  \\ \vdots \\ \sum_{i = 2}^\infty \frac{1}{i !} \left[ \overset{i}{\otimes}  \Delta x^T   \right] {\rm vec}( f_N^{(i)})     \end{bmatrix} 
\end{align}
where $\otimes$ is the Kronecker product, ${\rm vec}(\cdot)$ is an operator that stacks each column of a matrix into a single column vector, and higher order partial derivatives are defined recursively as
\begin{equation} \label{fkeq}
f_j^{(k)} = \frac{\partial {\rm vec} \big( f^{(k-1)}_j \big)}{\partial x^T} \added{ \in \mathbb{R}^{ N^{(k-1)} \times N}}.
\end{equation}
Above, all partial derivatives are evaluated at $x_0$ and the notation $ \overset{3}{\otimes}  \Delta x^T$  is shorthand for $\Delta x^T \otimes \Delta x^T \otimes \Delta x^T \added{\in \mathbb{R}^{1 \times N^3}}$.  By equating \eqref{timeexpand} and \eqref{taylfn}, substituting $\Delta x$ from \eqref{taylG} into \eqref{taylfn}, and matching the resulting powers in the isostable coefficients, one finds the resulting relationship
\begin{equation} \label{geqs}
0 = (J - ( \lambda_i + \lambda_j + \lambda_k +\dots ) {\rm Id} ) g^{ijk\dots} + q^{ijk\dots},
\end{equation} 
where ${\rm Id}$ is an appropriately sized identity matrix and $q^{ijk\dots} \in \mathbb{R}^N$ is a function of the lower order terms of the expansion, for example, $q^{321}$ (a third order term) might depend on $g^1$ or $g^{21}$ (a first or second order term) but will not depend on $g^{111}$ (a third order term).  Equation \eqref{geqs} is nearly identical to Equation (18) from \cite{wils20highacc} except for the fact that the time derivative of $g^{ijk\dots}$ is zero since the attractor is a fixed point.  For this reason,  $g^{ijk\dots}$ can be computed as the solution to a linear matrix equation

The individual terms of $I_n$ from \eqref{taylI} can be computed in a similar fashion.  Starting with equation (39) from \cite{wils20highacc}, for any isostable-based coordinate system the following holds for any state in the basin of attraction of the associated attractor:
\begin{equation}\label{prepsiadj}
\frac{d {I}_n}{dt}  = - \left(  \left. \frac{\partial F}{\partial x}^T\right|_x - \lambda_n {\rm Id}  \right) {I}_n. 
\end{equation}
Direct differentiation of \eqref{taylI} provides the time derivative of the gradient of the isostable coordinate
\begin{align}\label{myderivs}
\frac{d I_n}{dt} &=  \sum_{k = 1}^M  \left[\psi_k  \lambda_k  I_n^k \right]  + \sum_{j = 1}^M \sum_{k = 1}^j  \left[ \psi_j \psi_k  (\lambda_j+\lambda_k)  I_n^{jk}    \right]
 + \dots .
\end{align}
Above, the relationship $\dot{\psi}_k = \lambda_k \psi_k$ is used when taking time derivatives using the product rule.  Additionally, letting $J^T \equiv \frac{\partial F}{\partial x}^T = \begin{bmatrix} \frac{\partial f_1}{\partial x}^T &  \cdots  & \frac{\partial f_N}{\partial x}^T \end{bmatrix}$, where all partial derivatives are evaluated at $x_0$,  near the fixed point one can write the asymptotic expansion 
\begin{align} \label{dfdxexp}
\left.   \frac{\partial F}{\partial x}^T \right|_{x_0+\Delta x} &= J^T + \begin{bmatrix} a_1 & \cdots & a_N  \end{bmatrix}, \nonumber \\
 a_i &= \sum_{j = 1}^\infty \frac{1}{j !} \left( \left[ \overset{j}{\otimes}  \Delta x^T  \right] \otimes  {\rm Id}   \right){\rm vec}( f_i^{(j+1)}),
\end{align}
where $a_i$ is a column vector and $f_i^{(i+1)}$ was defined in \eqref{fkeq}.  Following the strategy from \cite{wils20highacc}, one can substitute \eqref{dfdxexp}, \eqref{myderivs}, and \eqref{taylI} into \eqref{prepsiadj}, replace any remaining $\Delta x$ terms with $G(\psi_1,\dots,\psi_M)$ from \eqref{taylG}, and match the resulting terms of the isostable coefficients to find the relationships
\begin{equation} \label{isorelations}
0 = -(J^T+ (-\lambda_n +\lambda_i + \lambda_j + \lambda_k + \dots){\rm Id} )I_n^{ijk\dots} + q_{\psi_n}^{ijk\dots},
\end{equation}
where  $q_{\psi_n}^{ijk\dots}$, is comprised of only lower order terms of the expansion of ${I}_n(\psi_1,\dots,\psi_M)$, for example, $q_{\psi_n}^{321}$ (a third order term) might depend on $I_n^1$ or $I_n^{21}$ (a first or second order term) but will not depend on $I_n^{111}$ (a third order term).  Equation \eqref{isorelations} is nearly identical to Equation (44) from \cite{wils20highacc}, except for the fact that the time derivative of $I_n^{ijk\dots}$ is zero since the attractor is a fixed point.  Once again, each $I_n^{ijk\dots}$ can be found as the solution to a linear matrix equation.  

A detailed description of strategies for computation of the terms $q^{ijk\dots}$ and $q_{\psi_n}^{ijk\dots}$ from equations \eqref{geqs} and \eqref{isorelations}, respectively, is discussed in \cite{wils20highacc} and is not repeated here.

\renewcommand{\thetable}{B\arabic{table}}  
\renewcommand{\thefigure}{B\arabic{figure}} 
\renewcommand{\theequation}{B\arabic{equation}} 
\setcounter{equation}{0}
\setcounter{figure}{0}

\section{Steady State Behavior of the Isostable Coordinate Expansion In Response to Sinusoidal Input} \label{apxb}

The relationships \eqref{ord1sol} and \eqref{ord2sol} illustrate how  the steady state behavior of the first and second order accurate terms of the isostable coordinate expansions in response to sinusoidal input $u(t) = \epsilon \sin(\omega t)$ can be written as a Fourier series with a fundamental frequency $\omega$ using a finite number of terms.  Theorem B.1~below shows that this pattern holds for all orders of accuracy in the expansion.

\vskip .1 in
\noindent {\bf Theorem B.1}:  Under the application of $u(t) = \epsilon \sin(\omega t)$, in steady state, any $\psi_{n,ss}^{(j)}$ term from \eqref{psiexp} can be written in the following form,
\begin{align} \label{psiss}
\psi_{n,ss}^{(j)} =  \sum_{k = 0}^{j} \bigg[ \alpha_{n,j}^k(\omega) \sin(k \omega t) + \beta_{n,j}^k(\omega) \cos( k \omega t) \bigg],
\end{align}
where $ \alpha_{n,j}^k(\omega) \in \mathbb{C}$ and $\beta_{n,j}^k(\omega) \in \mathbb{C}$.  

\vskip .1 in
\noindent {\it Proof}. Theorem B.1 can be proven by induction.  From equations \eqref{ord1sol} and \eqref{ord2sol}, this relationship holds for each $\psi_{n,ss}^{(1)}$ and $\psi_{n,ss}^{(2)}$, respectively.  Suppose that \eqref{psiss} holds for each $\psi_{n,ss}^{(1)},\psi_{n,ss}^{(2)}, \psi_{n,ss}^{(j-1)}, \dots, \psi_{n,ss}^{(j)}$ and  for all $n$.  Then considering the dynamics specified by \eqref{reducedmodinput} and using the asymptotic expansions \eqref{tinput} and \eqref{psiexp}, for any $\psi_n$ one can collect all $O(\epsilon^{j+1})$ terms giving the steady state dynamics
\begin{align} \label{sumterms}
\dot{\psi}_{n,ss}^{(j+1)} &= \lambda_n \psi_{n,ss}^{(j+1)} + \sum_{b_1 = 1}^M \bigg( \sin(\omega t) \psi_{b_1,ss}^{(j)} I_{n,A}^{b_1} \bigg) + \sum_{a_1 + a_2 = j}  \bigg(  \sin(\omega t) \bigg[  \sum_{b_1 = 1}^M \sum_{b_2 = 1}^{b_1}  \big[  \psi_{b_1,ss}^{(a_1)} \psi_{b_2,ss} ^{(a_2)}  I_{n,A}^{b_1 b_2}  \big] \bigg] \bigg)+ \dots  \nonumber \\
& +  \sum_{a_1 +a_2 + \dots + a_j = j} \bigg( \sin(\omega t) \bigg[   \sum_{b_1 = 1}^M \sum_{b_2 = 1}^{b_1} \cdots  \sum_{b_{j-1} = 1}^{b_{j-2}}    \sum_{b_j = 1}^{b_{j-1}}   \big[  \psi_{b_1,ss}^{(a_1)}  \dots \psi_{b_j,ss}^{(a_j)} I_{n,A}^{b_1\dots b_j}  \big]    \bigg] \bigg).
\end{align}
Consider, any term of the form $\sin(\omega t) \psi_{b_1,ss}^{(a_1)}  \dots \psi_{b_m,ss}^{(a_m)}$ for which $a_1 +a_2 + \dots + a_m = j$ from \eqref{sumterms}.  From \eqref{psiss}, terms of this form can be written as the product of sines and cosines with frequencies that are multiplies of $\omega$.  Considering the trigonometric product-to-sum identities
\begin{align} \label{prodid}
\sin(u) \sin(v) &= \frac{1}{2}( \cos(u-v) - \cos(u+v) ), \nonumber \\
\cos(u)\cos(v) &= \frac{1}{2}( \cos(u-v) + \cos(u+v) ), \nonumber \\
\sin(u)\cos(v) &= \frac{1}{2}( \sin(u+v) + \sin(u-v) ).
\end{align}  
By converting all products to sums, and noting, for instance, that $\psi_{b_1,ss}^{(a_1)}$ terms  have a maximum frequency of $a_1 \omega$, one can write
\begin{equation} \label{sindecomp}
\sin(\omega t) \psi_{b_1,ss}^{(a_1)}  \dots \psi_{b_m,ss}^{(a_m)} =  \sum_{k = 0}^{j+1} \bigg[ \alpha_k(\omega) \sin(k \omega t) + \beta_k(\omega) \cos( k \omega t) \bigg],
\end{equation}
where $\alpha_k(\omega) \in \mathbb{C}$ and $\beta_k(\omega) \in \mathbb{C}$.  Using the product-to-sum formulas, all terms containing sines and cosines in \eqref{sumterms} can be decomposed into the form \eqref{sindecomp}.  Therefore, \eqref{sumterms} can be rewritten in the form
\begin{equation} \label{psinss}
\dot{\psi}_{n,ss}^{(j+1)} = \lambda_n \psi_{n,ss}^{(j+1)}  + \sum_{k = 1}^{j+1} \bigg[  \gamma_{n,j+1}^k(\omega) \sin(k \omega t) + \delta_{n,j+1}^k (\omega) \cos(k \omega t)  \bigg],
\end{equation}
where each $ \gamma_{n,j}^k (\omega) \in \mathbb{C}$ and $\delta_{n,j}^k (\omega) \in \mathbb{C}$.  Finally, noticing that \eqref{psinss} is simply a periodically forced linear ordinary differential equation, one finds that the steady state dynamics are
\begin{equation} \label{steadyeq}
\psi_{n,ss}^{(j+1)}(t) = \sum_{k = 1}^{j+1} \frac{- \gamma_{n,j+1}^k(\omega) (k \omega  \cos(k \omega t) + \lambda_n \sin(k \omega t))  - \delta_{n,j+1}^k (\omega) (\lambda_n \cos(k \omega t) - k \omega \sin(k \omega t))   }{\lambda_n^2 + (k \omega)^2},
\end{equation}
which is of the form \eqref{psiss} thereby completing the proof. 
\vskip .1 in
\noindent {\it Remark}.  Theorem B.1~says that in response to an input $u(t) = \epsilon \sin(\omega t)$, that the $O(\epsilon^j)$ terms can be represented exactly as a Fourier series with a fundamental frequency $\omega$ and a maximum frequency $j \omega$.  This can be exploited to identify the terms of the expansions from   \eqref{taylgobs} and \eqref{tinput} as explained in the following section.

\renewcommand{\thetable}{C\arabic{table}}  
\renewcommand{\thefigure}{C\arabic{figure}} 
\renewcommand{\theequation}{C\arabic{equation}} 
\setcounter{equation}{0}
\setcounter{figure}{0}

\section{Fitting Isostable Reduced Models to Arbitrary Order Accuracy using Input-Output Relationships} \label{apxc}

Linear relationships \eqref{forelate} and \eqref{sorelate} can be used to fit the first and second order terms of the expansion \eqref{taylgobs} and \eqref{tinput}, respectively, when considering a system with a single output.  Similar relationships can be obtained for any desired orders of accuracy.  To show how this can be done, first consider the following Lemma:
\vskip .1 in
\noindent {{\bf Lemma C.1}}:  Under the application of $u(t) = \epsilon \sin(\omega t)$, for a single output system, in steady state the output relationship from \eqref{reducedoutput} can be written in the following form
\begin{align} \label{c1lemma}
y_{ss}(\omega,t) - y_0 &= \sum_{k = 0}^{j}  \bigg[  \tau_k (\omega) \sin(k \omega t) + \sigma_k (\omega) \cos(k \omega t)  \bigg]  + O(\epsilon^{j+1}),
\end{align}
where each $\tau_{j}(\omega) \in \mathbb{C}$ and $\sigma_{j} (\omega) \in \mathbb{C}$) is an $O(\epsilon^m)$ term with $m \geq j$.



 \vskip .1 in
\noindent {\it Remark}.  Lemma C.1 highlights that under the application of $u(t) = \epsilon \sin(\omega t)$, oscillatory components of frequency $j \omega$ of the steady state output are $O(\epsilon^m)$ terms where $m \geq j$. 

\vskip .1 in
\noindent {\it Proof}.  Consider any term from the Taylor expansion of $G_y(\psi_1,\dots,\psi_M)$ from \eqref{taylgobs} of the form $   \sum_{b_1 = 1}^M \sum_{b_2 = 1}^{b_1} \cdots  \sum_{b_{p-1} = 1}^{b_{p-2}}    \sum_{b_p = 1}^{b_{p-1}}   \big[  \psi_{b_1} \dots \psi_{b_p}  g_y^{b_1\dots b_p}  \big]$, where $p \in \mathbb{N}$.  Considering the asymptotic expansions of the isostable coordinates from Equation \eqref{psiexp}, in steady state one can rewrite these sums as
\begin{align} \label{expandlemma}
   \sum_{b_1 = 1}^M \sum_{b_2 = 1}^{b_1} \cdots  \sum_{b_{p-1} = 1}^{b_{p-2}} &   \sum_{b_p = 1}^{b_{p-1}}   \big[  \psi_{b_1,ss} \dots \psi_{b_p,ss}  g_y^{b_1\dots b_p}  \big]   =  \nonumber \\
   & \sum_{k = 1}^{j-1} \bigg(  \sum_{a_1 + a_2 + \dots + a_p=k} \bigg[  \epsilon^k      \sum_{b_1 = 1}^M \sum_{b_2 = 1}^{b_1} \cdots  \sum_{b_{p-1} = 1}^{b_{p-2}}    \sum_{b_p = 1}^{b_{p-1}}   \big[  \psi_{b_1,ss}^{(a_1)} \dots \psi_{b_p,ss}^{(a_p)}  g_y^{b_1\dots b_p}  \big]       \bigg] \bigg)  \nonumber \\
   &  +    \sum_{a_1 + a_2 + \dots + a_p=j} \bigg[  \epsilon^j      \sum_{b_1 = 1}^M \sum_{b_2 = 1}^{b_1} \cdots  \sum_{b_{p-1} = 1}^{b_{p-2}}    \sum_{b_p = 1}^{b_{p-1}}   \big[  \psi_{b_1,ss}^{(a_1)} \dots \psi_{b_p,ss}^{(a_p)}  g_y^{b_1\dots b_p}  \big]       \bigg]   + O(\epsilon^{j+1})
\end{align}
Above, all terms have been rewritten to emphasize the orders of $\epsilon$ in the Taylor expansion.  Using \eqref{psiss}, in steady state,  terms of the form $\psi_{n,ss}^{(j)}$ can be written as a Fourier series expansion with fundamental frequency $\omega$ and maximum frequency of $j \omega$.  Considering the product-to-sum identities from \eqref{prodid}, this implies Equation \eqref{expandlemma} can be rewritten as
\begin{align} \label{lemma1endpre}
  \sum_{b_1 = 1}^M  \sum_{b_2 = 1}^{b_1} & \cdots   \sum_{b_{p-1} = 1}^{b_{p-2}}    \sum_{b_p = 1}^{b_{p-1}}   \big[  \psi_{b_1,ss} \dots \psi_{b_p,ss}  g_y^{b_1\dots b_p}  \big]   =  \nonumber \\
 &  \sum_{k = 0}^{j-1}  \bigg[  \tau_{k,1}(\omega) \sin(k \omega t) + \sigma_{k,1}(\omega) \cos(k \omega t)  \bigg] \nonumber \\
  &  +    \sum_{a_1 + a_2 + \dots + a_p=j} \bigg[  \epsilon^j      \sum_{b_1 = 1}^M \sum_{b_2 = 1}^{b_1} \cdots  \sum_{b_{p-1} = 1}^{b_{p-2}}    \sum_{b_p = 1}^{b_{p-1}}   \big[  \psi_{b_1,ss}^{(a_1)} \dots \psi_{b_p,ss}^{(a_p)}  g_y^{b_1\dots b_p}  \big]       \bigg]   + O(\epsilon^{j+1}),
 \end{align}
 where each $\tau_{k,1}(\omega)$ and $\sigma_{k,1}(\omega) \in \mathbb{C}$ contain the contributions from the $O(\epsilon),\dots,O(\epsilon^{j-1})$ terms of the expansion \eqref{expandlemma}.  To arrive at the simplification \eqref{lemma1endpre} one can note that when using the product-to-sum identities \eqref{prodid} when simplifying each $\mathcal{O}(\epsilon^k)$ term for $k<j$, $k \omega$ is the maximum possible frequency of the resulting sinusoids.  Considering the $O(\epsilon^j)$ terms in \eqref{lemma1endpre},  continuing the simplification in this manner yields
 \begin{align}\label{lemma1end}
   \sum_{b_1 = 1}^M & \sum_{b_2 = 1}^{b_1}  \cdots   \sum_{b_{p-1} = 1}^{b_{p-2}}    \sum_{b_p = 1}^{b_{p-1}}   \big[  \psi_{b_1,ss} \dots \psi_{b_p,ss}  g_y^{b_1\dots b_p}  \big]   =  \nonumber \\
 &= \sum_{k = 0}^{j-1}  \bigg[  \tau_{k,2}(\omega) \sin(k \omega t) + \sigma_{k,2}(\omega) \cos(k \omega t)  \bigg] +  \bigg[ \tau_{j,2}(\omega) \sin(j \omega t) + \sigma_{j,2}(\omega) \cos(j \omega t) \bigg] + O(\epsilon^{j+1}).
\end{align}
where  terms that are proportional to $\sin(k \omega t)$ (resp.~$\cos(k \omega t)$) are absorbed into  the new constants $\tau_{k,2}(\omega)  \in \mathbb{C}$  (resp.~$\sigma_{k,2}(\omega) \in \mathbb{C})$.  Above, note that  $\tau_{j,2}(\omega)$ and $\sigma_{j,2}(\omega)$ are $O(\epsilon^j)$ terms.  Because any of the terms of the expansion \eqref{taylgobs} can be written in the form \eqref{lemma1end} it is possible to write $y_{ss}(\omega,t) - y_0$ in the form of \eqref{c1lemma} where the terms proportional to $\sin(j \omega t)$ and $\cos(j \omega t)$ are $O(\epsilon^j)$ terms which completes the proof.


\vskip .1 in
Lemma C.1 will used to show how higher order terms of the expansions \eqref{taylgobs} and \eqref{tinput} can be determined from the steady state behavior.  Towards this goal, \eqref{psiss} will be rewritten as
\begin{equation} \label{psirewrite}
\psi_{n,ss}^{(j)}(t) =\alpha_{n,j}^j(\omega) \sin(j \omega t) + \beta_{n,j}^j(\omega)  \cos( j \omega t)   +   \sum_{k = 0}^{j-1} \bigg[ \alpha_{n,j}^k(\omega)  \sin(k \omega t) + \beta_{n,j}^k(\omega)  \cos( k \omega t) \bigg],
\end{equation}
in order to emphasize the highest frequency terms in each $\psi_{n,ss}^{(j)}$.  The the following theorem  can now be stated:

\vskip .1 in
\noindent {{\bf Theorem C.1}}:   Consider an isostable coordinate reduction of the form \eqref{reducedmodinput} with output equation \eqref{reducedoutput} with a single output and $M$ isostable coordinates where all $\lambda_n$ for $n = 1 \dots M$ are known.  Suppose for all $k = 1,\dots,j$ that the coefficients $\alpha_{n,k}^k(\omega) $ and $\beta_{n,k}^k(\omega) $ from \eqref{psirewrite} have been computed; these represent the coefficients of \eqref{psirewrite} that are proportional to the highest frequency terms of each $\psi_{n,ss}^{(j)}(t)$.  Suppose also that all terms $g_y^{b_1}, g_y^{b_1 b_2}, \dots, g_y^{b_1 \dots b_j}$  and  $I_n^0, I_n^{b_1} ,\dots, I_n^{b_1,\dots,b_{j-1}}$ have been estimated for all $n$, i.e., all terms associated with the $j^{\rm th}$ and ${j-1}^{\rm th}$ order accurate expansions in the isostable coordinates from \eqref{taylgobs} and \eqref{tinput}, respectively, have been estimated.  Then
\begin{enumerate}[{Statement C}1:]
\item By applying a series of $q$ inputs of the form $u = \epsilon \sin(\omega t)$ with $w = w_1, \dots, w_q$ and measuring the resulting steady state behavior $y_{ss}(\omega,t)$, all terms of the form $g_y^{b_1 \dots b_{j+1}}$ and $I_n^{b_1,\dots,b_{j}}$ for all $n$ can be estimated using the linear relationship
\begin{equation} \label{arblin}
 \Upsilon_{j+1} =  \Xi_{j+1}^\dagger \bigg( \frac{1}{ \epsilon^{j+1}\pi}  \Gamma_{j+1} - R_{j+1}   \bigg)  + O(\epsilon),
\end{equation}
where 
\begin{equation} \label{gammaj1}
\Gamma_{j+1} =\begin{bmatrix}
\omega_1 \int_0^{2\pi/\omega_1} y_{ss}(w_1,t) \sin((j+1)\omega_1 t) dt   \\ \omega_1 \int_0^{2\pi/\omega_1} y_{ss}(w_1,t) \cos((j+1)\omega_1 t) dt  \\ \vdots \\
\omega_q \int_0^{2\pi/\omega_q} y_{ss}(w_q,t) \sin((j+1)\omega_q t) dt   \\ \omega_q \int_0^{2\pi/\omega_q} y_{ss}(w_q,t) \cos((j+1)\omega_q t) dt
\end{bmatrix} \in \mathbb{C}^{2q},
\end{equation}
where $\Upsilon_{j+1}\in \mathbb{C}^{\eta_{j+1}}$ is a vector containing $\eta_{j+1}$ total unknown elements and is comprised solely of terms of the form $g^{b_1 \dots b_{j+1}}$ and $I_n^{b_1,\dots,b_{j}}$, and $\Xi_{j+1} \in \mathbb{C}^{2q  \times \eta_{j+1}}$, $R_{j+1} \in \mathbb{C}^{2q}$ are appropriately defined matrices comprised solely of previously estimated terms, and $^\dagger$ denotes the pseudoinverse.
\item The terms $\alpha_{n,j+1}^{j+1}(\omega)$ and $\beta_{n,j+1}^{j+1}(\omega)$ of  $\psi_{n,ss}^{(j+1)}$ from the expansion \eqref{psirewrite} can be estimated for all $n$ once $\Upsilon_{j+1}$ has been estimated according to \eqref{arblin}.
\end{enumerate}

\vskip .1 in
\noindent {\it Proof of Statement C1.} 
Consider the steady state dynamics of $\psi_{n,ss}^{(j+1)}$ governed by \eqref{sumterms}.  These dynamics can be rewritten as a sum of sines and cosines according to \eqref{psinss}.  
\begin{align} \label{psinss_apxc}
\dot{\psi}_{n,ss}^{(j+1)} &= \lambda_n \psi_{n,ss}^{j+1}  +  \gamma_{n,j+1}^{j+1}(\omega) \sin(({j+1}) \omega t) + \delta_{n,j+1}^{j+1}(\omega) \cos(({j+1}) \omega t)  \nonumber \\
&+ \sum_{k = 1}^{j} \bigg[  \gamma_{n,j+1}^k(\omega) \sin(k \omega t) + \delta_{n,j+1}^k(\omega) \cos(k \omega t)  \bigg],
\end{align}
in the analysis to follow, the main focus will be on the terms proportional to $\sin((j+1) \omega t)$ and $\cos((j+1) \omega t)$ in the above equation.  When substituting \eqref{psirewrite} into \eqref{sumterms} and subsequently  applying the product-to-sum identities from \eqref{prodid}, the only way to yield terms that are proportional to $\sin((j+1) \omega t)$ and $\cos((j+1) \omega t)$ is to consider the highest frequency terms of each $\psi_{n,ss}^{(k)}$ (i.e.,~the highest frequency terms only depend on $\alpha_{n,k}^k(\omega)$ and $\beta_{n,k}^k(\omega)$ for $k = 1,\dots,j$).   With this in mind, considering the structure of \eqref{sumterms} and comparing to the structure of \eqref{psinss_apxc}, one finds that the terms $\gamma_{n,j+1}^{j+1}$ and $\delta_{n,j+1}^{j+1}$ can be written as
\begin{align} \label{highestcoeff}
\gamma_{n,j+1}^{j+1}(\omega) & = \sum_{b_1 = 1}^M \big[ \chi_{n,s}(\omega) I_{n,A}^{b_1} \big] +  \sum_{b_1 = 1}^M \sum_{b_2 = 1}^{b_1}  \big[ \chi_{n,s}^{b_1 b_2}(\omega) I_{n,A}^{b_1 b_2 } \big] +   \dots  \nonumber \\
&+  \sum_{b_1 = 1}^M \sum_{b_2 = 1}^{b_1} \cdots  \sum_{b_{j-1} = 1}^{b_{j-2}}    \sum_{b_j = 1}^{b_{j-1}}  \big[ \chi_{n,s}^{b_1 \dots b_j} (\omega) I_{n,A}^{b_1 \dots b_j} \big], \nonumber \\
\delta_{n,j+1}^{j+1}(\omega) & = \sum_{b_1 = 1}^M \big[ \chi_{n,c}^{b_1}(\omega) I_{n,A}^{b_1} \big] + \sum_{b_1 = 1}^M \sum_{b_2 = 1}^{b_1}  \big[ \chi_{n,c}^{b_1 b_2}(\omega) I_{n,A}^{b_1 b_2 } \big] +  \dots  \nonumber \\
&+  \sum_{b_1 = 1}^M \sum_{b_2 = 1}^{b_1} \cdots  \sum_{b_{j-1} = 1}^{b_{j-2}}    \sum_{b_j = 1}^{b_{j-1}}   \big[ \chi_{n,c}^{b_1 \dots b_j} (\omega)  I_{n,A}^{b_1 \dots b_j} \big],
\end{align}
where each term of the form $\chi_{n,s}^{b_1 b_2\dots}(\omega) \in \mathbb{C}$ and $\chi_{n,c}^{b_1 b_2\dots}(\omega) \in \mathbb{C}$ only depend on the highest frequency terms of each $\psi_{n,ss}^{(k)}$ (i.e.,~only depend on $\alpha_{n,k}^k(\omega)$ and $\beta_{n,k}^k(\omega)$ for $k = 1,\dots,j$ which are known by assumption).    Because  \eqref{psinss_apxc} is of the form \eqref{psinss} (i.e.,~a linear, periodically forced ordinary differential equation) the solution ${\psi}_{n,ss}^{(j+1)}(t) $ is given by \eqref{steadyeq}.   Note that from $\gamma_{n,j+1}^{j+1}(\omega)$ and $\delta_{n,j+1}^{j+1}(\omega)$ in \eqref{highestcoeff}, the only unknown terms are of the form $I_{n,A}^{b_1 \dots b_j}$.    Substituting \eqref{highestcoeff} into \eqref{steadyeq}, and recalling that $\alpha_{n,j+1}^{j+1}(\omega)$ and $\beta_{n,j+1}^{j+1}(\omega)$ are the coefficients proportional to the highest frequency terms of $\psi_{n,ss}^{(j+1)}$, one can write
\begin{align} \label{alphajj}
\alpha_{n,j+1}^{j+1}(\omega) & = \sum_{b_1 = 1}^M \sum_{b_2 = 1}^{b_1} \cdots  \sum_{b_{j-1} = 1}^{b_{j-2}}    \sum_{b_j = 1}^{b_{j-1}}  \big[ z_{n,\alpha}^{b_1 \dots b_j}(\omega)  I_{n,A}^{b_1 \dots b_j} \big] + k_{n,\alpha}^{j+1}(\omega) , \nonumber \\
\beta_{n,j+1}^{j+1} (\omega)& = \sum_{b_1 = 1}^M \sum_{b_2 = 1}^{b_1} \cdots  \sum_{b_{j-1} = 1}^{b_{j-2}}    \sum_{b_j = 1}^{b_{j-1}}  \big[ z_{n,\beta}^{b_1 \dots b_j} (\omega) I_{n,A}^{b_1 \dots b_j} \big] + k_{n,\beta}^{j+1}(\omega) ,
\end{align}
where each term of the form $z_{n,\alpha}^{b_1 b_2\dots}(\omega)  \in \mathbb{C}$, $z_{n,\beta}^{b_1 b_2\dots} (\omega)  \in \mathbb{C}$, $k_{n,\alpha}^{j+1} (\omega)  \in \mathbb{C}$, and $k_{n,\beta}^{j+1}(\omega)   \in \mathbb{C}$ can be determined from the terms of the reduction that are already known by assumption.

Next, invoking Lemma C.1 and considering terms up to and including $O(\epsilon^{j+1})$, multiplying both sides of \eqref{c1lemma} by either $\sin((j+1) \omega t)$ or $\cos((j+1) \omega t)$ and integrating over one period yields
\begin{align}
\begin{bmatrix} \label{stat2}
\int_0^{2 \pi/\omega} y_{ss}(\omega t) \sin((j+1)\omega t) dt \\
\int_0^{2 \pi/\omega} y_{ss}(\omega t) \cos((j+1)\omega t) dt
\end{bmatrix}  = \frac{\pi   }{\omega} \begin{bmatrix}   \tau_{j+1}(\omega)  \\ \sigma_{j+1}(\omega)    \end{bmatrix} + O(\epsilon^{j+2}),
\end{align}
where $\tau_{j+1}(\omega)$ and $\sigma_{j+1}(\omega)$ are $O(\epsilon^{j+1})$ terms.  Focusing on the terms $\tau_{j+1}(\omega)$ and $\sigma_{j+1}(\omega)$,  letting $\epsilon^j (y_{ss}(\omega,t)^{(j)})$ denote the $O(\epsilon^j)$ terms of the steady state output, one can collect all $O(\epsilon^{j+1})$ terms from the output using the asymptotic expansion \eqref{taylgobs} to yield
\begin{align} \label{sumtermsout}
y_{ss}^{(j+1)}(\omega,t) &= y_0 + \sum_{b_1 = 1}^M  \psi_{b_1,ss}^{(j+1)} g_y^{b_1} + \sum_{a_1 + a_2 = j+1}   \bigg[  \sum_{b_1 = 1}^M \sum_{b_2 = 1}^{b_1}  \big[  \psi_{b_1,ss}^{(a_1)} \psi_{b_2,ss} ^{(a_2)}  g_y^{b_1 b_2}  \big] \bigg] + \dots  \nonumber \\
& + \sum_{a_1 +a_2 + \dots +a_{j+1} = j+1}   \bigg[   \sum_{b_1 = 1}^M \sum_{b_2 = 1}^{b_1} \cdots  \sum_{b_{j} = 1}^{b_{j-1}}    \sum_{b_{j+1} = 1}^{b_{j}}   \big[  \psi_{b_1,ss}^{(a_1)}  \dots \psi_{b_{j+1},ss}^{(a_{j+1})} g_y^{b_1\dots b_{j+1}}  \big]    \bigg].
\end{align}
By converting all products to sums using \eqref{prodid} and recalling from \eqref{psiss} that each $  \psi_{b_1,ss}^{(a_1)}$ has terms with a maximum frequency of $a_1 \omega$, one can write
\begin{align} \label{xssfin}
y_{ss}^{(j+1)}(\omega,t)  &=\bigg[ \tau_{j+1}(\omega) \sin((j+1) \omega t) \nonumber + \sigma_{j+1}(\omega) \cos((j+1) \omega t) \bigg] \nonumber \\
&+ \sum_{k = 0}^j\bigg[r_{k,s}(\omega) \sin(k \omega t) + r_{k,c}(\omega) \cos(k \omega t)  \bigg],
\end{align}
where $r_{k,c}(\omega) \in \mathbb{C}$ and $r_{k,s}(\omega) \in \mathbb{C}$ are appropriately defined $O(\epsilon^{j+1})$ constants that are proportional to the lower frequency terms.  Note that above, the $\tau_{j+1}(\omega)$ and $\sigma_{j+1}(\omega)$ are the same as the constants from \eqref{c1lemma}, i.e.,~they are the $O(\epsilon^{j+1})$ terms that are proportional to $\sin((j+1)\omega t)$ and $\sin((j+1)\omega t)$, respectively.

Once again, note that the terms proportional to $\sin((j+1) \omega t)$ and $\cos((j+1) \omega t)$ in \eqref{xssfin} when applying the product-to-sum identities to \eqref{sumtermsout} only result when considering the highest frequency terms of each $\psi_{n,ss}^{(k)}$.   With this in mind, terms that are proportional to $\sin((j + 1)\omega t)$ and $\cos((j + 1)\omega t)$ can be computed solely with knowledge of each $\alpha_{n,k}^k(\omega)$ and $\beta_{n,k}^k(\omega)$ for $k \leq j+1$.    By assumption, each $\alpha_{n,k}^k(\omega)$ and $\beta_{n,k}^k(\omega)$ is already known for  $k \leq j$ for all $n$.  Considering the relationship \eqref{alphajj}, the terms  $\alpha_{n,j+1}^{j+1}(\omega)$ and $\beta_{n,j+1}^{j+1}(\omega)$ are known for all $n$ except for the contributions from terms of the form $ I_{n,A}^{b_1 \dots b_j}$.

By substituting \eqref{psirewrite} into \eqref{sumtermsout}, expanding using product-to-sum identities \eqref{prodid}, and keeping in mind the structure of the unknown terms of $\alpha^{j+1}_{n,j+1}(\omega)$ and $\beta_{n,j+1}^{j+1}(\omega)$ from \eqref{alphajj}, one can collect terms appropriately as in \eqref{xssfin} to find that $\tau_{j+1}(\omega)$ and $\sigma_{j+1}(\omega)$ can be written as

\begin{align} \label{tauandsigma}
\tau_{j+1}(\omega) &= \sum_{n = 1}^M  \bigg[          \sum_{b_1 = 1}^M \sum_{b_2 = 1}^{b_1} \cdots  \sum_{b_{j-1} = 1}^{b_{j-2}}    \sum_{b_j = 1}^{b_{j-1}}  \big[ z_{n,I,\tau}^{b_1 \dots b_j}(\omega) I_{n,A}^{b_1 \dots b_j} \big]             \bigg]   \nonumber \\
&+     \sum_{b_1 = 1}^M \sum_{b_2 = 1}^{b_1} \cdots  \sum_{b_{j} = 1}^{b_{j-1}}    \sum_{b_{j+1} = 1}^{b_{j}}   \bigg[  z_{g,\tau}^{b_1\dots b_{j+1}}(\omega) g_y^{b_1\dots b_{j+1}}  \bigg]  + r_{\tau,j+1}(\omega), \nonumber \\
\sigma_{j+1}(\omega) &= \sum_{n = 1}^M  \bigg[          \sum_{b_1 = 1}^M \sum_{b_2 = 1}^{b_1} \cdots  \sum_{b_{j-1} = 1}^{b_{j-2}}    \sum_{b_j = 1}^{b_{j-1}}  \big[ z_{n,I,\sigma}^{b_1 \dots b_j}(\omega) I_{n,A}^{b_1 \dots b_j} \big]             \bigg] \nonumber \\
& +    \sum_{b_1 = 1}^M \sum_{b_2 = 1}^{b_1} \cdots  \sum_{b_{j} = 1}^{b_{j-1}}    \sum_{b_{j+1} = 1}^{b_{j}}  \bigg[  z_{g,\sigma}^{b_1\dots b_{j+1}}(\omega) g_y^{b_1\dots b_{j+1}}  \bigg]  + r_{\sigma,j+1}(\omega),
\end{align} 
where each $z_{n,I,\tau}^{b_1 b_2 \dots}(\omega) $, $z_{n,I,\sigma}^{b_1 b_2 \dots}(\omega)$, $z_{n,g,\tau}^{b_1 b_2 \dots}(\omega)$,  $z_{n,g,\sigma}^{b_1 b_2 \dots}(\omega)$, $r_{\tau,j+1}(\omega)$ and $r_{\sigma,j+1} (\omega)\in \mathbb{C}$ and can be computed from the lower order terms of the expansions \eqref{taylgobs} and \eqref{tinput} that are known by assumption.  Here, the terms $ r_{\tau,j+1}(\omega)$ and $r_{\sigma,j+1}(\omega)$ are neither proportional to the terms $g_y^{b_1\dots b_{j+1}}$ nor $g_y^{b_1\dots b_{j+1}}$.

By considering the steady state behavior in response to inputs of the form $u(t) = \epsilon \sin(\omega t)$ at multiple frequencies, one can substitute Equation \eqref{tauandsigma} into Equation  \eqref{stat2} and rearrange to yield
\begin{equation} \label{largeraccuracy}
\Gamma_{j+1} = \pi \epsilon^{j+1} \big( \Xi_{j+1} \Upsilon_{j+1} + R_{j+1}   \big) + O(\epsilon^{j+2}),
\end{equation}
where $\Gamma_{j+1}$ is defined in \eqref{gammaj1}, $\Upsilon_{j+1}\in \mathbb{C}^{\eta_{j+1}}$ is a vector containing $\eta_{j+1}$ total elements of the unknown terms $g_y^{b_1 \dots b_{j+1}}$ and $I_n^{b_1,\dots,b_{j}}$, $\Xi_{j+1}$ is a known matrix that contains terms that are proportional to unknown terms from $\Upsilon_{j+1}$,  and $R_{j+1}$ is a known vector comprised of all remaining terms.    The relationship \eqref{arblin} follows immediately.

\vskip .1 in
\noindent {\it Proof of Statement C2.}  
\noindent Estimation of  $\Upsilon_{j+1}$ according to \eqref{arblin} provides and estimation of all terms of the form $ I_{n,A}^{b_1 \dots b_j}$.  Consequently, all elements of the right hand side of \eqref{alphajj} are known providing an estimate for the terms $\alpha_{n,j+1}^{j+1}(\omega)$ and $\beta_{n,j+1}^{j+1}(\omega)$
 for all $n$.  

\vskip .1 in
\noindent {\it Remark.}
\noindent Theorem C.1~provides a strategy to iteratively estimate the terms of the expansions \eqref{taylgobs} and \eqref{tinput} that are part of the reduced order model \eqref{reducedmodinput} with output \eqref{reducedoutput}.  For example, equation \eqref{foestimate} provides a relatively straightforward expression of this form that can be used to infer the the first order accurate terms.   Equations for identifying second order and higher accuracy terms are more complicated and it is advisable to compute the terms of matrices $\Xi_{j}$ and $R_j$ using a symbolic computational package.  All equations of the form \eqref{arblin} require knowledge of all lower order terms of the expansions \eqref{taylgobs} and \eqref{tinput}.



\end{appendices}

\FloatBarrier


\end{document}